%% file: orbifold_dt.tex
\documentclass[12pt]{amsart}

\usepackage{times}

\usepackage{amsmath}
\usepackage{amssymb}
\usepackage{color}
\usepackage{epsfig}
\usepackage{graphicx}
\usepackage{subfig}
\usepackage{mathtools} 

\usepackage{setspace}

\newtheorem{theorem}{Theorem}[section]
\newtheorem{conjecture}[theorem]{Conjecture}
\newtheorem{definition}[theorem]{Definition}
\newtheorem{corollary}[theorem]{Corollary}
\newtheorem{lemma}[theorem]{Lemma}

\newcommand{\bbC}{\mathbb{C}}
\newcommand{\bbZ}{\mathbb{Z}}

\newtheorem{remark}[theorem]{Remark}

\newtheorem{proposition}[theorem]{Proposition}

\newcommand{\ZtwoZtwo}{\mathbb{Z}_{2}\times \mathbb{Z}_{2}}
\newcommand{\Hilb}{\operatorname{Hilb}}
\newcommand{\cnums} {{\mathbb C}}          
\newcommand{\znums} {{\mathbb Z}}		
\renewcommand{\O}{\mathcal{O}}
\renewcommand{\hat}{\widehat}

\newcommand{\vdim}{\operatorname{vdim}}

\begin{document}

\title[Coloured 3D Young diagrams]{Generating functions for coluored 3D Young diagrams and the Donaldson-Thomas invariants of orbifolds}
\author{Benjamin Young, with an appendix by Jim Bryan}
\date{\today}
\maketitle

\bibliographystyle{amsplain}

\begin{abstract}
We derive two multivariate generating functions for three-dimensional Young diagrams (also called plane partitions).  The variables correspond to a colouring of the boxes according to a finite abelian subgroup $G$ of $SO(3)$.  These generating functions turn out to be orbifold Donaldson--Thomas partition functions for the orbifold $[\bbC^3 / G]$.
We need only the vertex operator methods of Okounkov--Reshetikhin--Vafa for the easy case $G = \bbZ_n$; to handle the considerably more difficult case $G=\bbZ_2\times\bbZ_2$, we will also use a refinement of the author's recent $q$--enumeration of pyramid partitions.

In the appendix,  we relate the diagram generating
functions to the Don\-ald\-son-Thom\-as partition functions of the orbifold
$[\bbC^3 / G]$. We find a relationship between the Don\-ald\-son-Thom\-as
partition functions of the orbifold and its $G$-Hilbert scheme resolution. We
formulate a crepant resolution conjecture for the Don\-ald\-son-Thom\-as theory
of local orbifolds satisfying the Hard Lefschetz condition.

\end{abstract}

\onehalfspacing

\section{Introduction}

A \emph{3D Young diagram}, or \emph{3D diagram} for short, is a stable pile of cubical boxes which sit in the corner of a large cubical room.  More formally, a 3D Young diagram is a finite subset $\pi$ of $(\bbZ_{\geq 0})^3$ such that if any of 
\[(i+1,j,k), (i,j+1,k), (i,j,k+1)\]
are in $\pi$, then $(i,j,k) \in \pi$.
The ordered triples are the ``boxes''; the closure condition means that the boxes of a 3D partition are stacked stably in the positive octant, with gravity pulling them in the direction $(-1,-1,-1)$.  

3D Young diagrams are well--studied; they are also called \emph{plane partitions} or \emph{3D partitions} elsewhere in the literature.  The first result on 3D Young diagrams is due to Dr.~Percy MacMahon~\cite{MACMAHON}.  MacMahon was the first to ``$q$--count'' (i.e. to give a generating function for) 3D Young diagrams by volume:
\begin{equation}
\label{eqn:macmahon_gf}
\sum_{\pi \text{ 3D diagram}} q^{|\pi|} = \prod_n \left( \frac{1}{1-q^n}\right)^n,
\end{equation}
where $|\pi|$ denotes the number of boxes in $\pi$.  Generating functions of this form will appear frequently, so we adopt the following notation:
\begin{definition}
Let 
\begin{align*}
M(x,q) &= \prod_{n=1}^{\infty} \left(\frac{1}{1-xq^n}\right)^n \\
\widetilde{M}(x,q) &= M(x,q) M(x^{-1}, q)
\end{align*}
\end{definition}
We call $M(x,q)$ and $\widetilde{M}(x,q)$ the \emph{MacMahon} and \emph{MacMahon tilde} functions, respectively.  Strictly speaking, $\widetilde{M}(x,q)$ lies in the ring of formal power series $\bbZ[\![x, x^{-1}, q]\!]$.  However, in all of our applications, we will specialize $x$ and $q$ in such a way that no negative powers of any variables appear in the formulae (see Theorems~\ref{thm:zn_result}~and~\ref{thm:z2z2_result}).

Since MacMahon, there have been many proofs of (\ref{eqn:macmahon_gf}), spanning many fields: combinatorics, statistical mechanics, representation theory, and others.  Recently, there has been a thorough study of the various symmetry classes of 3D Young diagrams~\cite{BRESSOUD}, and of many macroscopic properties of large random 3D Young diagrams~\cite{Ok-Re}. There is also active research in algebraic geometry which relies upon enumerations of various types of 3D partitions~\cite{MNOP}.

We will derive two refinements of MacMahon's generating function. Fix a set of colours $\mathcal{C}$, and replace the variable $q$ with a set of variables,
\[
\mathcal{Q} = \{q_g\;|\;g \in \mathcal{C}\}.
\]
We will need to assign a colour to each point of the first orthant.  In particular, we will usually have $\mathcal{C} = G$, a finite Abelian group.  In this case, addition in $\bbZ_{\geq 0}^3$ must respect the group law of $G$.  

\begin{definition}
A \emph{colouring} is a map 
\[
K:(\bbZ_{\geq 0})^3 \rightarrow \mathcal{C}.
\]
If $\mathcal{C} = G$ is a finite Abelian group, then a \emph{$G$--colouring} is a colouring which is also a homomorphism of additive monoids.
\end{definition}
Note that a $G$--colouring is uniquely determined by $K(1,0,0)$, $K(0,1,0)$ and $K(0,0,1)$, and that $K(0,0,0)$ is the identity element of $G$.  

There is a simple way of defining a $G$--colouring $K_G$ when $G$ is a three--dimensional matrix group $G$.  Decompose $G$ as a direct sum of one--dimensional representations $R_x, R_y, R_z$.  
The set of irreducible representations of any Abelian $G$ forms a group $\widehat{G} \simeq G$ under tensor product, so let $\psi$ be an isomorphism $\psi:\widehat{G} \longrightarrow G $ and define
\[
K_G(i,j,k) = \psi(
R_x^{\otimes i} \otimes R_y^{\otimes j} \otimes R_z^{\otimes k}).
\]
Both of the colourings used in this paper are of this form.

We next define the multivariate generating function $Z_G = Z_G(\mathcal{Q})$ which ``$\mathcal{Q}$--counts'' diagrams (that is, $Z_G$ counts each diagram with the $\mathcal{Q}$--weight of its boxes):

\begin{definition}
For $g \in G$, let $|\pi|_g$ be the number of $g$--coloured boxes in $\pi$,
\[
|\pi|_g = |K_G^{-1}(g) \cap \pi|.
\]
Define the \emph{$G$-coloured partition function}   
\[
Z_G = \sum_{\pi \text{3D partition}} \prod_{g \in G} q_g^{|\pi|_g}.
\]
\end{definition}

The question of determining $Z_G$, though completely combinatorial,
has its genesis in a field of enumerative algebraic geometry called
\emph{Don\-ald\-son-Thom\-as theory}.  When $G$ is a finite Abelian subgroup
of $SO(3)$ (which forces $G = \bbZ_n$ or $\bbZ_2 \times \bbZ_2$),
there is a colouring induced by the natural three dimensional
representation for which the generating function $Z_G$ is, up to signs
of the variables, the orbifold Don\-ald\-son--Thom\-as partition function
for the quotient stack $[\bbC^3/G]$ (see Appendix~\ref{app: DT theory
of C3/G}).  Although it is not yet clear why, these seem to be
precisely the groups $G$ for which $Z_G$ has a product formula.

\begin{theorem}
\label{thm:zn_result}
Let $G = \bbZ_n$ and let the colouring $K_{Z_n}$ be given by
\begin{align*}
K_{Z_n}(1,0,0) &= 1 \\
K_{Z_n}(0,1,0) &= -1 \\
K_{Z_n}(0,0,1) &= 0.
\end{align*}
Let $q=q_0\cdots q_{n-1}$ and for $a,b \in [1, n-1]$, let $q_{[a,b]} = q_aq_{a+1} \cdots q_b$.
Then \[
Z_{\bbZ_n} = M(1, q)^n \prod_{0 < a \leq b < n} 
\widetilde{M}(q_{[a,b]}, q).
\]
\end{theorem}

The proof of Theorem~\ref{thm:zn_result} is straightforward; it is essentially a simple modification of the methods used in~\cite{ORV} (or, indeed, a special case of the extremely general methods of \cite{Ok-Re}).  We include it for completeness and as an introduction to the vertex operator calculus used to prove Theorem~\ref{thm:z2z2_result}.  There are several other ways to prove Theorem~\ref{thm:zn_result}, some of which have (at least implicitly) appeared in the literature.  For example, ~\cite{ANDREWS-PAULE, GANSNER} both compute a generating function with variables $x_k (k \in \bbZ)$ which can be easily specialized to $Z_{\bbZ_n}$.  The result~\cite{ANDREWS-PAULE} is particularly notable, as it is a direct computer algebra implementation of MacMahon's techniques of combinatory analysis.  The following theorem, however, is new:

\begin{theorem}
\label{thm:z2z2_result}
Let $G = \bbZ_2 \times \bbZ_2 = \{0, a, b, c\}$ and let the colouring $K_{\bbZ_2 \times \bbZ_2}$ be given by
\begin{align*}
K_{Z_2 \times Z_2}(1,0,0) &= a \\
K_{Z_2 \times Z_2}(0,1,0) &= b \\
K_{Z_2 \times Z_2}(0,0,1) &= c.
\end{align*}
Let $q=q_0q_aq_bq_c$.  Then
\[
Z_{\bbZ_2 \times \bbZ_2} = M(1,q)^4 \cdot \frac{
\widetilde{M}(q_aq_b, q)
\widetilde{M}(q_aq_c, q)
\widetilde{M}(q_bq_c, q)
}{
\widetilde{M}(-q_a, q)
\widetilde{M}(-q_b, q)
\widetilde{M}(-q_c, q)
\widetilde{M}(-q_aq_bq_c, q)
}.
\]
\end{theorem}

\begin{figure}
\caption{A partition coloured according to $K_{\bbZ_2 \times \bbZ_2}$ and to $K_{\bbZ_3}$}
\label{fig:colourings}
\begin{center}
\subfloat[$K_{\bbZ_2\times \bbZ_2}$ -- weight $q_0^{30}q_a^{29}q_b^{31}q_c^{28}$]{
	\includegraphics[width=2in]{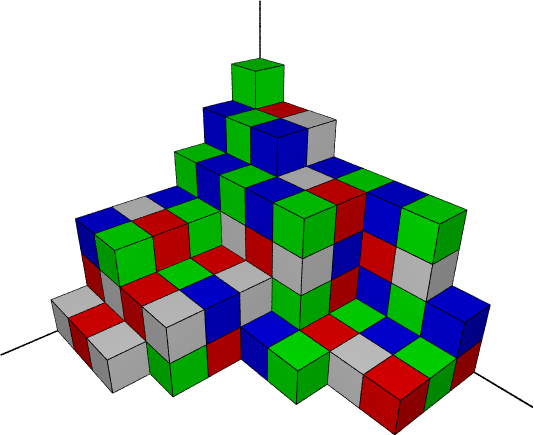}
}
\subfloat[$K_{\bbZ_3}$ -- weight $q_0^{40}q_1^{38}q_2^{40}$]{
	\includegraphics[width=2in]{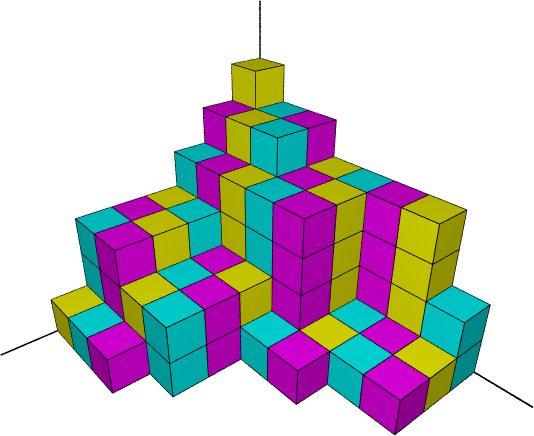}
}
\end{center}
\end{figure}

See Figure~\ref{fig:colourings} for pictures of a partition coloured in the manner described by these theorems.

As an application of these theorems, we will compute the Don\-ald\-son-Thom\-as invariants of the orbifolds $[\bbC^3/\bbZ_n]$ and $[\bbC^3/\bbZ_2 \times \bbZ_2]$. The orbifold Don\-ald\-son-Thom\-as partition function of $[\bbC^3/G]$ has variables labeled  by representations of $G$ (see Appendix) and hence has the same variables as the $G$-coloured diagram partition function. In the Appendix, we prove that the diagram partition function and the Don\-ald\-son-Thom\-as partition function are related by simple sign changes on the variables:
\begin{theorem}
\label{thm:orbifold_dt_partition_functions}
The orbifold Don\-ald\-son-Thom\-as partition functions of the orbifolds $[\bbC^3/\bbZ_2 \times \bbZ_2]$ and $[\bbC^3/\bbZ_n]$ are given by
\begin{align*}
Z^{\text{DT}}_{\bbC^3/\bbZ_n}(q_0,q_1,...,q_{n-1}) 
&= Z_{\bbZ_n}(-q_0,q_1,...,q_{n-1}) \\
Z^{\text{DT}}_{\bbC^3/\bbZ_2 \times \bbZ_2}(q_0,q_a,q_b,q_c) 
&= Z_{\bbZ_2 \times \bbZ_2}(q_0,-q_a,-q_b,-q_c) \\
\end{align*}
where $q$ and $q_{[a,b]}$ are defined as in Theorems~\ref{thm:zn_result} and~\ref{thm:z2z2_result}.
\end{theorem}
There is a striking similarity between the Don\-ald\-son-Thom\-as partition functions of the orbifold $[\bbC^3/G]$ and the crepant resolution given by the $G$--Hilbert scheme.  The following is proved in the Appendix:
\begin{theorem}
\label{thm:resolution_dt_partition_functions} Let $Y_G \longrightarrow
\bbC^3/G$ be the crepant resolution of $\bbC^3/G$ given by the
$G$-Hilbert scheme.  $Y_G$ has a natural basis of curve classes
indexed by non-trivial elements of $G$.  The Don\-ald\-son-Thom\-as
partition functions of $Y_{\bbZ_n}$ and $Y_{\bbZ_2 \times \bbZ_2}$ are
given by
\begin{align*}
Z^{\text{DT}}_{Y_{\bbZ_n}} &= M(1,-q)^n \prod_{0 < a \leq b < n} M(q_{[a,b]}, -q),  \\
Z^{\text{DT}}_{Y_{\bbZ_2 \times \bbZ_2}} &=
M(1,-q)^4
\frac{
M(q_aq_b, -q) M(q_bq_c, -q) M(q_aq_c, -q)
}{
M(q_a,-q)M(q_b,-q)M(q_c,-q)M(q_aq_bq_c,-q)
},
\end{align*}
where $\{q_1,...,q_{n-1}\}$ and $\{q_a,q_b,q_c\}$ are the variables corresponding to curve classes and $q$ is the variable corresponding to Euler number.
\end{theorem}
We see from these theorems that the reduced partition function of the
orbifold $[\bbC^3/G]$ is obtained from the reduced partition function
of the resolution by identifying the variables appropriately and then
simply writing a tilde over every factor of $M$ in the formula!
A similar phenomenon was observed by Szendr\H{o}i for the partition
function of the (non--commutative) conifold singularity and its
crepant resolution~\cite{NC_DT}.

It would be very desirable to have even a conjectural understanding of
the relationship between the Don\-ald\-son--Thom\-as theory of an arbitrary
Calabi--Yau orbifold and its crepant resolution(s).  We formulate a
conjecture for the case of a local orbifold satisfying the hard
Lefschetz condition (see Conjecture~\ref{conj: local DT CRC}.

Theorem~\ref{thm:z2z2_result} is not straightforward to prove.  Essentially none of the standard proofs of MacMahon's colourless result can be modified to work in this situation. 
The generating function was first conjectured by Jim Bryan based on some related phenomena from Don\-ald\-son--Thom\-as theory; concurrently,
Kenyon made an (unpublished) equivalent conjecture for $\bbZ_2 \times \bbZ_2$--weighted dimer models on the hexagon lattice, based on computational evidence.  

Having this conjectured formula was crucial for finding the proof of Theorem~\ref{thm:z2z2_result}, which involves a somewhat bizarre detour: one must first $\mathcal{Q}$--count \emph{pyramid partitions} (see Figure~\ref{fig:pyramid}). One then performs a computation with vertex operators to make $Z_{\bbZ_2 \times \bbZ_2}$ emerge.  We discovered this idea serendipitously while trying to generalize our earlier work on pyramid partitions ~\cite{MYSELF}.

\section{Review: the infinite wedge space} 

Our general strategy will be to think of a 3D diagram $\pi$ as a set of
diagonal slices, $\{\pi_k \;| \; k \in \bbZ\}$, where $\pi_k$ is the set of all
bricks which lie in the plane $x-y=k$.  We will then analyze how one passes from one slice to the next.  Since we will be summing over all 3D Young diagrams, it is very helpful to consider (possibly infinite) formal sums of the form
\[ 
\sum_{\lambda \in \text{some set of partitions} } f_{\lambda}(\mathcal{Q}) \cdot \lambda,
\]
where $f_{\lambda}(\mathcal{Q})$ is a power series in the elements of $\mathcal{Q}$.  A nice way of describing the set of all such sums is the \emph{charge--zero subspace of the infinite wedge space}, 
\[
(\Lambda^{\infty/2})_0 V
\]
where $V$ is a vector space with a basis labeled by the elements of $\bbZ+\frac{1}{2}$.  This setting allows one to define, quite naturally, several useful operators on partitions.  

\begin{figure}
\caption{Applying $\alpha_{-3}$ to a partition}
\label{fig:add_ribbon}
\begin{center}
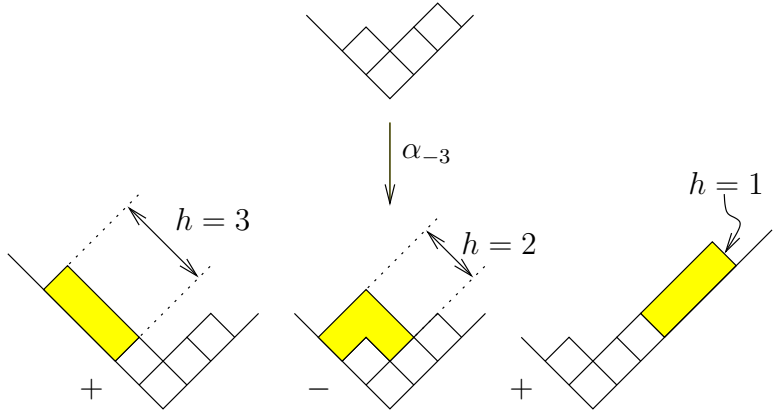 
\end{center}
\end{figure}

The use of $(\Lambda^{\infty/2})_0 V$, and its associated operators, was in part popularized by~\cite[Appendix A]{Ok-random-partitions}, and we shall adhere to the notation established there.  In this section, we have collected the minimum number of formulae necessary for our purposes.  We will use Dirac's ``bra--ket'' notation 
\[
\left<\lambda\left|\mu\right>\right.
\]
to denote the inner product under which the partitions are orthonormal.  We will need
need the bosonic creation and annihilation operators $\alpha_n$, defined in ~\cite[Appendix A]{Ok-random-partitions} in the section on Bosons and Vertex Operators.  The operators $\alpha_n$ satisfy the Heisenberg commutation relations,
\begin{equation}
\label{eqn:alpha_commutator}
[\alpha_n, \alpha_{-m}] = n\delta_{m,n}.
\end{equation}
Concretely, $\alpha_{-n}$ acts on a 2D Young diagram $\lambda$ by adding a single border strip of length $n$ onto $\lambda$ in all possible ways, with sign $(-1)^{h+1}$, where $h$ is the height of the border strip (see Figure~\ref{fig:add_ribbon}).  The operator $\alpha_n$ is adjoint to $\alpha_{-n}$, and acts by deleting border strips.

Let $x_j (j \geq 1)$ be indeterminates; and define the homogeneous, elementary, and power sum symmetric functions as usual:
\begin{align*}
\sum_i h_i(x_1, x_2, \ldots) t^i &= \prod_i \frac{1}{1-x_it}\\ 
\sum_i e_i(x_1, x_2, \ldots) t^i &= \prod_i (1+x_it)\\
p_i(x_1, x_2, \ldots) &= \sum_{j \geq 1} x_j^i 
\end{align*}
For a comprehensive reference on symmetric functions, see \cite{ECII}. 
We next define the vertex operators $\Gamma_{\pm}$:
\begin{definition}
\label{def:vertex_operator}
\[
\Gamma_{\pm}(x_1, x_2, \ldots) = \exp \sum_{k} \frac{p_k}{k} \alpha_{\pm k}
\]
\end{definition}
The matrix coefficients (with respect to the orthonormal basis formed by the 2D Young diagrams) of the $\Gamma_{\pm}$ operators turn out ~\cite[A.15]{Ok-random-partitions} to be the skew Schur functions,
\[
\left<\lambda \left| \Gamma_{-}(x_1, x_2, \ldots) \right| \mu\right> = 
\left<\mu \left| \Gamma_{+}(x_1, x_2, \ldots) \right| \lambda\right> = s_{\lambda/\mu}(x_i).
\]
We will need the following well-known theorem from representation theory (see, for example,~\cite[Chapter 8]{FH}) to work with $\Gamma_{\pm}$ and other exponentiated operators.
\begin{theorem} (Campbell-Baker-Hausdorff)
If $A$ and $B$ are operators, then 
\[
\log(\exp(A) \exp(B)) = A + B + \frac{1}{2}[A,B] + 
\cdots,
\]
where the higher--order terms are multiples of nested commutators of $A$ and $B$.
\end{theorem} 
It is certainly possible to give more terms in the expansion, but we shall only need the following two corollaries.

\begin{corollary}  
\label{cor:CBH_II}
If $A$ and $B$ are commuting operators, then 
\[
\exp (A) \exp (B)=\exp(A + B) .
\]
\end{corollary}

\begin{corollary} 
\label{cor:CBH_I}
If $A$ and $B$ are operators such that $[A,B]$ is a central element, then we have
\[
\exp(A) \exp(B) = \exp([A,B]) \exp(B) \exp(A).
\]
\end{corollary}

\section{The operators $\Gamma(x)$, $\Gamma'(x)$, and $Q_g$}
Our next goal is to define precisely what it means for two diagonal slices $\lambda, \mu$ to sit next to one another in a 3D Young diagram, and to define operators for working with such slices.  

\begin{definition} Let $\lambda, \mu$ be two 2D Young diagrams.  We say that $\lambda$ \emph{interlaces with} $\mu$, and write $\lambda \succ \mu$, if $\mu \subseteq \lambda$ and the skew diagram $\lambda / \mu$ contains no vertical domino.
\end{definition}
For example, $(6,3,2) \succ (4,2)$, because the skew diagram $(6,3,2)/(4,2)$ has no two boxes in the same column.
The following lemma is easy to check:
\begin{lemma} The following are equivalent:
\begin{enumerate}
\label{lem:interlacing}
\item $\lambda \succ \mu$.
\item The row lengths ${\lambda_i}, {\mu_i}$ satisfy $\lambda_1 \geq \mu_1 \geq \lambda_2 \geq \mu_2 \geq \cdots$.
\item \label{sublemma:col_interlace} $\lambda^i - \mu^i = 0 \text{ or } 1$, for each pair of \emph{columns} $\lambda^i, \mu^i$.
\item \label{sublemma:why_interlace}$\lambda$ and $\mu$ are two adjacent diagonal slices of some 3D Young diagram.
\end{enumerate}
\end{lemma}
Note that we have used the convenient, but slightly nonstandard, notation $\lambda^i$ to denote the \emph{columns} of $\lambda$.

Part (\ref{sublemma:col_interlace}) will become relevant in Section~\ref{sec:pyramid_part}, when we will see that adjacent diagonal slices of a pyramid partition also interlace.

We are mainly interested in two specializations of $\Gamma_{\pm}(x_1,x_2,\ldots)$ which create interlacing partitions, and which depend only upon a single indeterminate $q$.  The first will be denoted $\Gamma_{\pm}(q)$, and is obtained by performing the specialization $x_1 \mapsto q$, $x_i \mapsto 0$ for $i>1$.  Its formula is 
\begin{equation}
\label{eqn:gamma_formula}
\Gamma_{\pm}(q) = \exp \sum_k \frac{q^k}{k} \alpha_{\pm k}.
\end{equation}
Recall ~\cite[Chapter 7]{ECII} that if $\lambda,\mu$ are partitions, then we may define the \emph{skew Schur function} $s_{\lambda/\mu}(x_1,x_2, \cdots)$ by $\sum_T x^T$,
where $T$ runs over the set of semi\-standard tableaux of shape $\lambda/\mu$.
Following~\cite{ORV}, we see that
\[
s_{\lambda/\mu}(q, 0, 0, \ldots) = 
\begin{cases}
q^{|\lambda|-|\mu|} & \text{if }\lambda \succ \mu \\
0 & \text{if }\lambda \not \succ \mu. \\
\end{cases}
\]
One can then show~\cite{Ok-Re-skew} that
\begin{align}
\label{eqn:gamma_interpretation}
\Gamma_-(q) \mu &= \sum_{\lambda \succ \mu}q^{|\lambda| - |\mu|} \lambda  &
\Gamma_+(q) \lambda &= \sum_{\mu \prec \lambda}q^{|\lambda| - |\mu|}\mu. 
\end{align}

For the second specialization, recall that there is an involution $\omega$ on the algebra of symmetric functions \cite[Chapter 7.6]{ECII}, given by any one of the following equivalent definitions:
\begin{align*}
e_k(x_i) &\longleftrightarrow h_k(x_i) \\
p_k(x_i) &\longleftrightarrow (-1)^{k-1}p_k(x_i) \\
s_{\lambda/\mu}(x_i) &\longleftrightarrow s_{\lambda'/\mu'}(x_i) 
\end{align*}
Here, $\lambda'$ is the transpose partition of $\lambda$.
To obtain the second specialization, called $\Gamma_{\pm}'(q)$, we first perform the involution $p_k \mapsto \omega p_k$ and then specialize $x_1 \mapsto q, x_i \mapsto 0 \;(i>1)$ as before.  We obtain the formula
\begin{equation}
\label{eqn:second_gamma_formula}
\Gamma_{\pm}'(q) = \exp \sum_k \frac{(-1)^{k-1}q^k}{k} \alpha_{\pm k}
\end{equation}
with the property that
\begin{align*}
\label{eqn:gamma_interpretation}
\Gamma_-'(q) \mu &= \sum_{\lambda' \succ \mu'}q^{|\lambda| - |\mu|} \lambda  &
\Gamma_+'(q) \lambda &= \sum_{\mu' \prec \lambda'}q^{|\lambda| - |\mu|}\mu. 
\end{align*}

\begin{lemma} If $a$ and $b$ are commuting variables, then we have the following multiplicative commutators in $\mathbb{C}[\![a,b]\!]$:
\begin{align*}
[\Gamma_{+}(a), \Gamma_{-}'(b)] &= {1+ab}  &
[\Gamma_{+}'(a), \Gamma_{-}(b)] &= {1+ab}  \\
[\Gamma_{+}(a), \Gamma_{-}(b)] &= \frac{1}{1-ab}  & 
[\Gamma_{+}'(a), \Gamma_{-}'(b)] &= \frac{1}{1-ab} \\ 
\end{align*}
\end{lemma}

\begin{proof}
Let us compute the first of these commutators; the others are similar.  Let us apply~(\ref{eqn:alpha_commutator}), and then use Corollary~\ref{cor:CBH_I} to rephrase the answer as the exponential of a commutator. We have 
\begin{align*}
[\Gamma_{+}'(a), \Gamma_{-}(b)] &=
\exp \sum_{j,k} \frac{(-1)^{j-1}a^jb^k}{jk} [\alpha_j, \alpha_{-k}] \\
&= \exp\left( -\sum_{j} \frac{(-ab)^j}{j} \right)\\
&= \exp(\log(1 - (-ab)) ).
\end{align*}
\end{proof}

We next define diagonal operators $Q_g$ for assigning weights to 2D partitions.
\begin{definition}
\label{def:weight_operator} For $g \in G$, define the \emph{weight operator} $Q_g$ by
\[
Q_g \left | \lambda\right > = q_g^{|\lambda|} \left| \lambda \right>.
\]
\end{definition}
The operator $Q_g$ can be commuted past any of the $\Gamma_{\pm}$ operators, at the expense of changing the argument of $\Gamma_{\pm}$:
\begin{align*}
\Gamma_+(x) Q_g &= Q_g \Gamma_+(xq_g) &
Q_g \Gamma_-(x) &= \Gamma_-(xq_g) Q_g    \\
\Gamma_+'(x)Q_g &= Q_g \Gamma_+'(xq_g) &
Q_g\Gamma_-'(x) &= \Gamma_-'(xq_g) Q_g  .
\end{align*}

\section{Counting with $Z_n$ colouring}
\label{sec:zn_count}
As a motivating example, let us use (\ref{eqn:gamma_interpretation}) to write down a vertex operator expression which computes MacMahon's generating function~(\ref{eqn:macmahon_gf}), using the variable $q=q_0$.  This formula appears in ~\cite{ORV} with marginally different notation.

Consider a 3D Young diagram $\pi$ and its diagonal slices:
\[\phi \prec \cdots \prec \pi_{-2} \prec \pi_{-1} \prec \pi_0 
\succ \pi_1 \succ \pi_2 \succ \cdots \succ \phi,
\]
where $\phi$ denotes the empty partition.
Each such $\pi$ contributes 
\[
q_0^{|\pi|}=q_0^{\sum |\pi_n|}
\] 
to the generating function, so we have
\[
\sum_{\pi \text{ 3D diagram}}q^{|\pi|} = \left< \phi \left|
\prod_{i=1}^{\infty}\left(
\Gamma_+(1) Q_0
\right)
\prod_{i=1}^{\infty}\left(
\Gamma_-(1) Q_0
\right)
\right| \phi \right>.
\]
This works because the operators $\Gamma_-$ and $\Gamma_+$ pass from one slice to the next larger (respectively smaller) slice in all possible ways, and the $Q_0$ operators assign the proper weight to each slice.  One then commutes all the $\Gamma_-$ operators to the left and all the $\Gamma_+$ operators to the right (following the method outlined in~\cite{ORV}) to compute the generating function.

\begin{figure}
\caption{Slicing a $\bbZ_3$--coloured 3D diagram}
\label{fig:zn_slice}
\begin{center}
\subfloat{
	\includegraphics[width=2in]{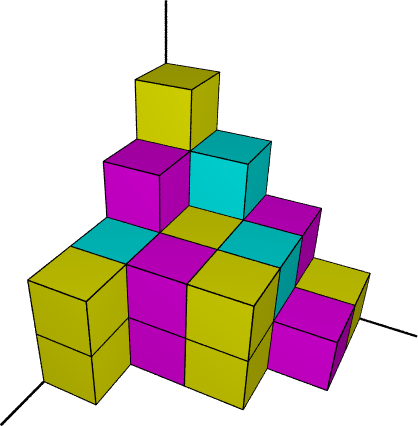}
}
\subfloat{
	\includegraphics[width=2in]{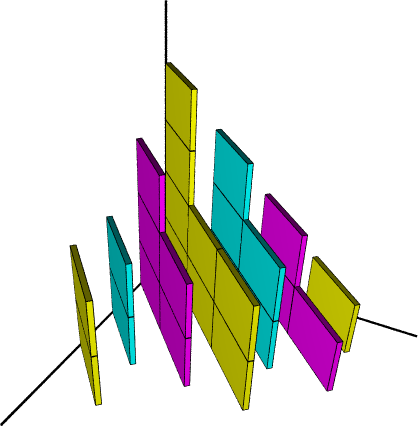}
}
\end{center}
\end{figure}

Let us now write down a vertex operator expression which computes $Z_{\bbZ_n}$.  Here, $\mathcal{Q} = \{q_0, \ldots, q_{n-1}\}$, $q=q_0q_1\cdots q_{n-1}$, and $K = K_{\bbZ_n}$.  The computation is straightforward (following precisely the method of~\cite{ORV}) but awkward, so it is helpful to organize the work by collecting together $n$ vertex operators at a time.  Note that the diagonal slices of $\pi$ are all monochrome (see Figure~\ref{fig:zn_slice}), so we define 
\[
\overline{A}_\pm(x) = 
\Gamma_{\pm}(x) Q_1
\Gamma_{\pm}(x) Q_2
\cdots
Q_{n-1}
\Gamma_{\pm}(x) Q_0
\]

Then, the following vertex operator product counts $\bbZ_{n}$--coloured 3D diagrams:
\begin{equation}
\label{eqn:zn_vertex_product}
Z_{\bbZ_n}  = \left< \phi \left| \cdots \overline{A}_+(1)\overline{A}_+(1)\overline{A}_+(1) \overline{A}_-(1)\overline{A}_-(1)\overline{A}_-(1) \cdots \right| \phi \right>
\end{equation}
Let $q = q_0q_1\cdots q_{n-1}$, and let $Q = Q_0Q_1\cdots Q_{n-1}$.  We use the commutation relations of the previous section to compute
\begin{align*}
\overline{A}_+(x) &= Q \cdot
\Gamma_+\left(xq_1q_2q_3\cdots q_{n-1}q_0\right)
\Gamma_+\left(xq_2q_3 \cdots q_{n-1}q_0\right) \cdots 
\Gamma_+\left(xq_0\right) \\
\overline{A}_-(x) &= 
\Gamma_+(x)
\Gamma_+(xq_1)
\cdots
\Gamma_+(xq_1q_2 \cdots q_{n-1})
\cdot Q\\
&=
\Gamma_+\left(xq q_1^{-1}q_2^{-1}\cdots q_{n-1}^{-1}q_0^{-1}\right)
\Gamma_+\left(xq q_2^{-1}q_3^{-1}\cdots q_{n-1}^{-1}q_0^{-1}\right) \cdots 
\\ & \quad \quad
\cdots
\Gamma_+\left(xqq_0^{-1}\right)
\cdot Q
\end{align*}
Next, set 
\begin{align*}
A_+(x) &= Q^{-1}\overline{A}_+(x); & A_-(x) &= \overline{A}_-(x)Q^{-1}.
\end{align*}
From this expression, it is clear that 
\[
A_+(x)A_-(y) = C(x,y) \cdot A_-(y)A_+(x)
\]
where $C(x,y)$ is the following product of the $n^2$ commutators obtained by moving a $\Gamma_+$ past a $\Gamma_-$:
\[
\begin{split}
C(x,y) = \left(
\frac{1}{1-qxy}
\right)^n  
\prod_{0 \leq a \leq b < n}
\left(
\frac{1}{1-(q_aq_{a+1}\cdots{q_b})qxy}
\right) \\ \cdot 
\prod_{0 \leq a \leq b < n}
\left(
\frac{1}{1-(q_aq_{a+1}\cdots{q_b})^{-1}qxy}
\right)
\end{split}
\]

  We now follow the derivation of MacMahon's formula in~\cite{ORV}.  Starting with (\ref{eqn:zn_vertex_product}), we convert all of the $\overline{A}_{\pm}$ into $A_{\pm}$ and move the resulting weight functions to the outside of the product (where they act trivially).  This gives
\[
Z_{\bbZ_n} = \left< \phi \left|
\cdots A_+(q^2) A_+(q) A_+(1) 
A_-(1) A_-(q) A_-(q^2) \cdots
\right|\phi\right>.
\]
We then commute all $A_+$ operators to the right and all $A_-$ to the left:
\begin{align*}
\label{eqn:Zn_vertex_ops}
Z_{\bbZ_n}
&=\left< \phi \right| \cdots A_+(q^2)A_+(q)\underbracket{A_+(1) A_-(1)}A_-(q)A_-(q^2) \cdots \left| \phi \right> \\
&=  C(1,1)\left< \phi \left| \cdots A_+(q^2)\underbracket{A_+(q)A_-(1) A_+(1)A_-(q)}A_-(q^2) \cdots \right| \phi \right> \\
&= \cdots \\
&=  \prod_{i,j = 0}^{\infty}\!\!C(q^i,q^j) \, \cdot \, \left< \phi \left| A_-(1)A_-(q)A_-(q^2) \cdots A_+(q^2) A_+(q)A_+(1) \right| \phi \right>. \\
\end{align*}
The vertex operator product in the final line is now equal to 1, because $\left<\phi \right|A_-(x) = \left< \phi \right|$ and $A_+(x)\left|\phi\right> = \left|\phi \right>$.  Finally, we rewrite the remaining product with MacMahon functions:
\begin{align*}
Z_{\bbZ_n} &= 
\prod_{i,j=0}^{\infty} C(q^i,q^j)  \\
&= M(1,q)^n \!\! \prod_{0 \leq a \leq b < n} M(q_a\cdots q_b, q) M(q_a^{-1}\cdots q_b^{-1}, q) \\
&= M(1,q)^n \!\! \prod_{0 \leq a \leq b < n} \widetilde{M}(q_a\cdots q_b, q),
\end{align*}
and Theorem~\ref{thm:zn_result} is proven. $\hfill \square$
\section{Pyramid partitions}
\label{sec:pyramid_part}

\begin{figure}
\caption[A pyramid partition, removed from the pyramid]{}
\label{fig:pyramid}
\begin{center}
\subfloat[The set $\mathcal{B}$ of bricks]{\includegraphics[width=2.4in]{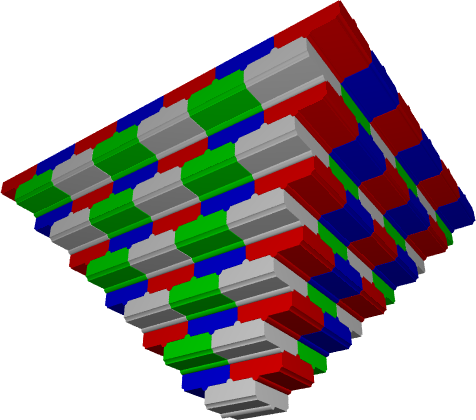}}
\subfloat[A pyramid partition removed from $\mathcal{B}$]{\includegraphics[width=2.4in]{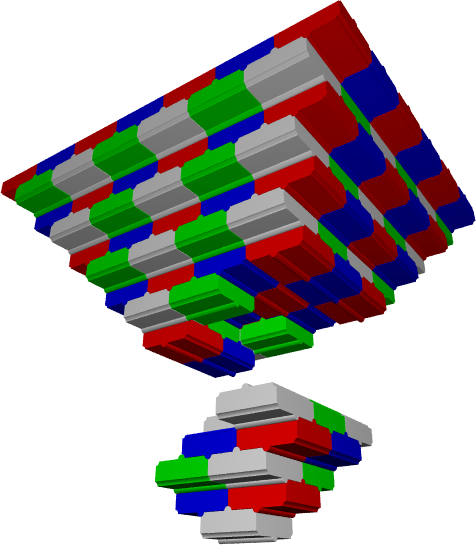}}
\end{center}
\end{figure}

The methods of the Section~\ref{sec:zn_count} may also be used to $\mathcal{Q}$-count a similar type of three-dimensional combinatorial object, called \emph{pyramid partitions}.   Essentially, we want to replace $\bbZ_{\geq 0}^3$ with the upside--down pyramid shaped stack of bricks shown in Figure~\ref{fig:pyramid}.   Note that the bricks have ridges and grooves set into them; this helps to remind us how the bricks are meant to stack.  

Szendr\H{o}i~\cite{NC_DT} introduced us to the ideas in this section, albeit in
a different context.  He proves that counting pyramid partitions with a
slightly simpler colour scheme (namely specializing $q_0 = q_c, q_b=q_a$) yields
a certain noncommutative Don\-ald\-son--Thom\-as partition function.   We shall
borrow some of Szendr\H{o}i's terminology, but not much of the machinery that
he developed.

\begin{figure}
\caption[A quiver diagram for colouring pyramid partitions]{}
\label{fig:pyramid_partition}
\begin{center}
\subfloat[The quiver $P$] { 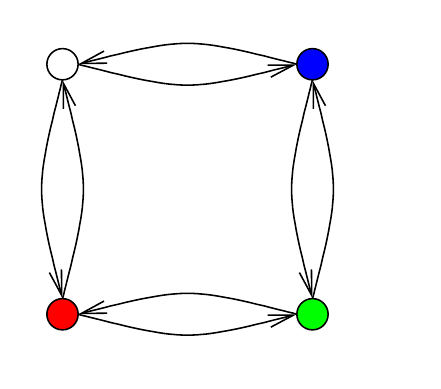 }
\subfloat[A pyramid partition $\pi$]{\includegraphics[width=2in]{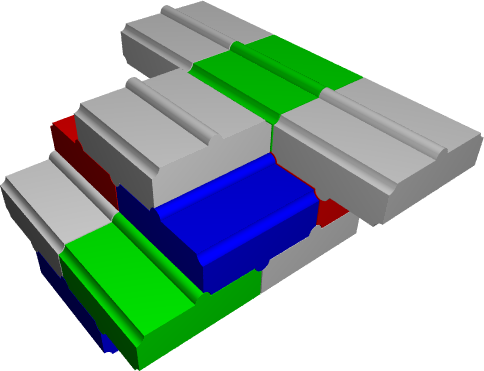}}
\end{center}
\end{figure}

We will start by giving a rather algebraic definition for the bricks in a pyramid partition.  Consider the quiver (or directed graph) $P$ shown in Figure~\ref{fig:pyramid_partition}(a).  The vertices of $P$ are the elements of $\bbZ_2 \times \bbZ_2 = \{0,a,b,c\}$.  The edges are labelled $\{v_1,w_1,v_2,w_2\}$. 

\begin{definition} A \emph{word} in $P$ is the concatenation of the edge labels of some directed path in $P$.  We may optionally associate a \emph{base} to a word; the base is the starting vertex of the path.
\end{definition}
Note that a word based at $0$ may also be based at $c$, but not at $b$ or $a$.
Any path in $P$ is uniquely determined by its base and its word.
\begin{definition} 
Form the \emph{path algebra} $\bbC P$ spanned by all words in $P$, and define the noncommutative quotient ring $A = \bbC P / I_W$, where
\begin{align*}
I_W &= \left<v_1w_iv_2 - v_2w_iv_1, w_1v_jw_2 - w_2v_jw_1 \right>,  &
i,j & \in \{1,2\}.
\end{align*}
If $B$ is a word in $\bbC P$, we write $[B]$ for its residue class in $\bbC P / I_w$.
\end{definition}

\begin{definition}
A \emph{brick} is an element $[B]$ of $\bbC P / I_W$, where $B$ is a word based at the vertex 0.  Let $\mathcal{B}$ be the set of all bricks.
\end{definition}

To understand how to draw Figure~\ref{fig:pyramid}, we interpret the edge labels of $P$ as vectors in $\bbZ^3$.
\begin{definition} Let
\begin{align*}
v_1 &= (-1,1,0) & v_2 &= (1,1,0)	\\
w_1 &= (0,1,-1) & w_2 &= (0,1,1).
\end{align*}
The \emph{position} of a brick $[B]$ is the sum of the vectors corresponding to the edge labels in $[B]$.  The brick corresponding to the empty walk, $[\:]$, is located at the origin.
\end{definition}

We next define a ``colouring'' on $\mathcal{B}$.
\begin{definition} Define
\[
K_{\text{pyramid}}:\mathcal{B} \longrightarrow \bbZ_2 \times \bbZ_2
\]
by setting $K_{\text{pyramid}}([B])$ to be the final vertex of any path whose word is $B$.  We call $K_{\text{pyramid}}([B])$ the \emph{colour} of $B$.
\end{definition}

For an example of all of these concepts, define the brick $[B]$ by the word $B = v_2w_2v_2w_1$.  The brick $[B]$ is based at the vertex $0$, ending at the vertex $b$.  The position of $[B]$ is $(2,4,0)$; $[B]$ is the $c$-coloured brick in the top layer of Figure~\ref{fig:pyramid_partition}b.  

\begin{figure}
\caption[A view of the pyramid $\mathcal{B}$ from below]{A view of the pyramid $\mathcal{B}$ from below.  The bricks have been shrunk to points, and some checkerboard--coloured slices are shown.}
\label{fig:pyramid_below}
\begin{center}
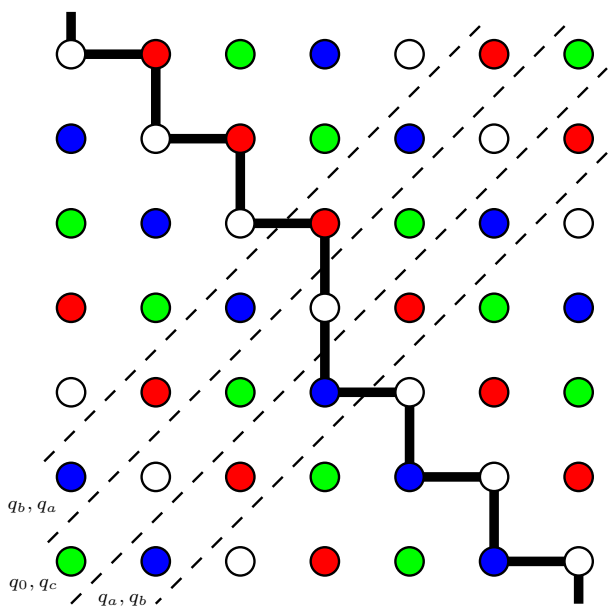
\end{center}
\end{figure}

  Note that the colour is completely determined by the $x$ and $y$ coordinates of $[B]$; Figure~\ref{fig:pyramid_below} shows the colouring as viewed from along the $z$ axis.

\begin{definition}
\label{def:pyramid_partition}
A \emph{pyramid partition} $\pi$ is a subset of $\mathcal{B}$ such that if $[B] \in \pi$ then every prefix of $B$ also represents a brick in $\pi$.  
\end{definition}

Note that pyramid partitions may also be defined algebraically, although it is unnecessary to do so for this paper. A pyramid partition corresponds to a \emph{framed cyclic $\bbC P / I_W$--module based at 0}, much in the same way that a 3D Young diagram corresponds to a monomial ideal in $\bbC[x,y,z]$.  We refer the reader to~\cite{NC_DT} for further details of this approach.

  Our next goal is to show that the diagonal slices of a pyramid partition interlace with one another.  This will allow us to reuse the strategy of Section~\ref{sec:zn_count} to obtain a nice generating function for pyramid partitions.

\begin{definition} Let $\pi$ be a pyramid partition.  Define the \emph{$k$th diagonal slice of $\pi$}, written $\pi_k$, to be the set of all bricks in $\pi$ whose position $(x,y,z)$ satisfies $x-y=k$.
\end{definition}

\begin{lemma} 
\label{lem:pyramid_slice}
Let $k \geq 0$.  Then
\begin{align*}
\pi_{-2k} &= \{ [(v_1w_2)^k W] \} \cap \pi,	   \\
\pi_{2k} &= \{ [(v_2w_1)^k W] \} \cap \pi, 	  
\end{align*}
where $W$ runs over all words in $v_1w_1$ and $v_2w_2$, and 
\begin{align*}
\pi_{-2k-1} &= \{ [(v_1w_2)^kv_1  W'] \} \cap \pi, \\
\pi_{2k+1} &= \{ [(v_2w_1)^kv_2 W'] \} \cap \pi, 
\end{align*}
where $W'$ runs over all words in $w_1v_1$ and $w_2v_2$.  Moreover, the bricks of $\pi_k$ form a 2D Young diagram; the slices are single--coloured, as follows:
\[
\text{Colour of }\pi_k = 
\begin{cases}
0 & \text{if }k = 0 \pmod{4} \\
b & \text{if }k = 1 \pmod{4} \\
c & \text{if }k = 2 \pmod{4} \\
a & \text{if }k = 3 \pmod{4}.
\end{cases}
\]
\end{lemma}
\begin{figure}
\caption[A diagonal slice of a pyramid partition is a 2D Young diagram]{A diagonal slice of a pyramid partition, interpreted as a 2D Young diagram}
\label{fig:pyramid_singleslice}
\begin{center}
	\includegraphics[height=2in]{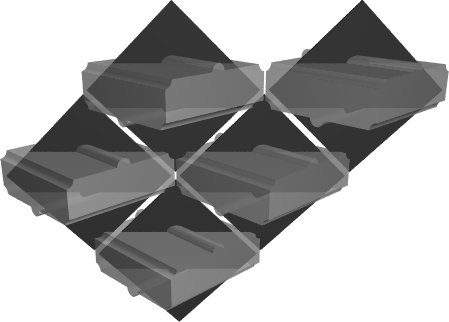}
\end{center}
\end{figure}

\begin{proof}
Let us prove the first equation; the other three are similar.  Note that the brick represented by the word $(v_1w_2)^k$ is in position $(-k,2k,k)$ and thus lies in the $-2k$th diagonal.  Appending $v_1w_1$ or $v_2w_2$ to this word adds $(-1,2,-1)$ or $(1,2,1)$ to the position, which does not alter $x-y$.  

To see that the bricks of $\pi_{-2k}$ form a 2D Young diagram, observe that $w_1v_1$ and $w_2v_2$ commute in $\bbC P /I_w$.  The suffix $(w_1v_1)^i(w_2v_2)^j$ corresponds to the $(i,j)$ box in the Young diagram.  Again, the other cases are similar.  The colours are easy to check directly. 
\end{proof}
Figure~\ref{fig:pyramid_singleslice} shows the central slice $\pi_0$ of the pyramid partition of Figure~\ref{fig:pyramid_partition}b.  Every brick has been replaced with a square tile to make the orientation of the 2D Young diagram clear. 

\begin{lemma} 
\label{lem:pyramid_interlacing}
For $k \geq 0$, we have the following interlacing properties (where the prime denotes transposition of 2D Young diagrams):
\begin{align*}
\pi_{2k} &\succ \pi_{2k+1}    &
\pi'_{-2k} &\succ \pi'_{-2k-1}   \\
\pi'_{2k+1} &\succ \pi'_{2k+2}  & 
\pi_{-2k-1} &\succ \pi_{-2k-2}  
\end{align*}
\end{lemma}

\begin{figure}
\caption[Row-- and column--interlacing behaviour in a pyramid partition]{Row-- and column--interlacing behaviour for adjacent diagonal slices of a pyramid partition}
\label{fig:pyramid_doubleslice}
\begin{center}
\subfloat{
	\includegraphics[width=2in]{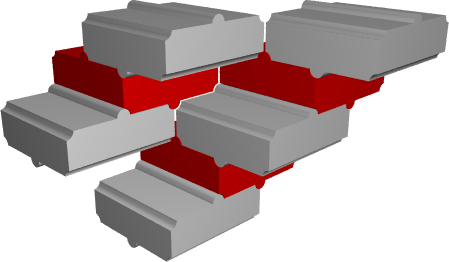}
}
\subfloat{
	\includegraphics[width=2in]{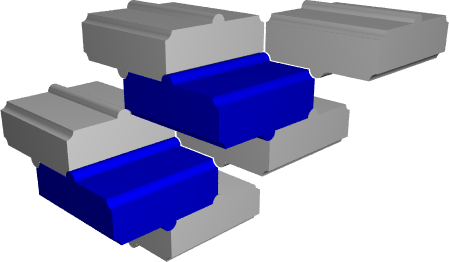}
}
\end{center}
\end{figure}

\begin{proof}
Let us handle the first case.  Let $R^j_{k}$ be the set of bricks in the $j$th column of $\pi_{k}$, and suppose that $|R^j_{k}| = \ell^j_k$.  Explicitly,
\begin{align*}
R^j_{2k+1} &= \{ (v_2w_1)^kv_2(w_1v_1)^j(w_2v_2)^i \;|\; 0 \leq i < \ell^j_{2k+1} \}\\
 &= \{ (v_2w_1)^k(v_1w_1)^j(v_2w_2)^iv_2 \;|\; 0 \leq i < \ell^j_{2k+1} \},\\
R^j_{2k} &= \{ (v_2w_1)^k(v_1w_1)^j(v_2w_2)^i \;|\; 0 \leq i < \ell^j_{2k} \}.\\
\end{align*}

In particular, each of the bricks in $R^j_{2j} \cup R^j_{2j+1}$  may be represented as some prefix of the word
\[
(v_2w_1)^k(v_1w_1)^j(v_2w_2)^{\max\{\ell^j_{2k},\ell^j_{2k+1}\}}.
\]

Informally speaking, $R^j_{2k}$ and $R^j_{2k+1}$ form a chain of bricks, each of which rests on the previous one (see Figure~\ref{fig:interlace}).  It follows from Definition~\ref{def:pyramid_partition} that $\ell^j_{2k} - \ell^j_{2k+1} \in \{0, 1\}$; then part (\ref{sublemma:col_interlace}) of Lemma~\ref{lem:interlacing} says that $\pi_{2k} \succ \pi_{2k+1}$.

\begin{figure}
\caption{The $j$th columns of two adjacent slices $\pi_0$, $\pi_1$}
\label{fig:interlace}
\begin{center}
	\includegraphics[height=2in]{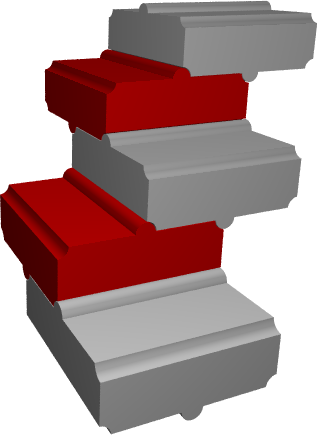}
\end{center}
\end{figure}

Next let us see that $\pi'_{2k+1} \succ \pi'_{2k+2}$.  Let $R_{i, k}$ be the $i$th \emph{row} of $\pi_k$, with $|R_{i, k}| = \ell_{i, k}$.  We have
\begin{align*}
R_{i,2k+2} &= \{ (v_2w_1)^k(v_2w_2)^i(v_1w_1)^j \;|\; 0 \leq j < \ell_{i, 2k+2} \}\\
 &= \{ (v_2w_1)^kv_2(w_2v_2)^i(w_1v_1)^j \;|\; 0 \leq j < \ell_{i, 2k+2} \},\\
R_{i, 2k+1} &= \{ (v_2w_1)^kv_2(w_2v_2)^i(w_1v_1)^j \;|\; 0 \leq j < \ell_{i,2k+1} \}.\\
\end{align*}
from which it follows that $\ell_{j,2k+1} - \ell_{j,2k+2} \in \{0,1\}$.  This means that $\pi'_{2k+1} \succ \pi'_{2k+2}$.  See Figure~
\ref{fig:pyramid_doubleslice} for an illustration of the difference between the row--interlacing and column--interlacing behaviour.  

The remaining cases are similar.
\end{proof}

\section{A generating function for pyramid partitions}

We will now compute the following generating function for pyramid partitions.
\begin{definition}
Let $\pi$ be a pyramid partition.  For $g \in \bbZ_2 \times \bbZ_2$, let
\[ 
|\pi|_g = |K_{\text{pyramid}}^{-1}(g) \cap \pi|
\]
denote the number of boxes coloured $g$ in $\pi$.  Define
\[
Z_{\text{pyramid}} = 
\sum_{\pi \text{ pyramid partition}} \prod_{g \in \bbZ_2 \times \bbZ_2} q_g^{|\pi|_g}.
\]
\end{definition}

\begin{theorem}
\label{thm:pyramid_gf}
\[
Z_{\text{pyramid}}
= 
\frac{
   M(1,q)^4 \widetilde{M}(q_bq_c)\widetilde{M}(q_aq_c)
}{
   \widetilde{M}(-q_a, q) \widetilde{M}(-q_b, q)
   \widetilde{M}(-q_c, q) \widetilde{M}(-q_aq_bq_c, q) 
},
\]
where $q = q_0q_aq_bq_c$.
\end{theorem}

Theorem~\ref{thm:pyramid_gf} may seem unrelated to the other theorems in this paper, but it will turn out that it is the key to computing $Z_{\bbZ_2 \times \bbZ_2}$.  The proof is much like that of Theorem~\ref{thm:zn_result}.

\begin{proof}
We define a vertex operator product which counts pyramid partitions.  Let us first define an operator which sweeps out four slices of the pyramid partition at the same time. Let
\[
\overline{A}'_{\pm}(x) = 
\Gamma_{\pm}(x)Q_b\Gamma_{\pm}'(x)Q_c
\Gamma_{\pm}(x)Q_a\Gamma_{\pm}'(x)Q_0,
\]
so that 
\[
Z_{\text{pyramid}} 
=
\left< \phi \left|
\cdots \overline{A}'_+(1) \overline{A}'_+(1) \overline{A}'_-(1) \overline{A}'_-(1) \cdots
\right| \phi \right>. 
\]
It is simple to check this product against Lemmas~\ref{lem:pyramid_interlacing} and~\ref{lem:pyramid_slice} to be sure that it describes the correct colouring and interlacing behaviour.  Set
\begin{align*}
A'_+(x) &= Q_0^{-1}Q_b^{-1}Q_c^{-1}Q_a^{-1}\overline{A}'_+(x) \\
A'_-(x) &= \overline{A}'_+(x)Q_0^{-1}Q_b^{-1}Q_c^{-1}Q_a^{-1}.
\end{align*}
Commuting the weight operators past the vertex operators gives
\[
A'_+(x) = \Gamma_+(xq_bq_cq_aq_0)
\Gamma_+'(xq_cq_aq_0)
\Gamma_+(xq_aq_0)
\Gamma_+'(xq_0) 
\]
and
\[
\begin{split}
A'_-(y) = 
\Gamma_-(yqq_b^{-1}q_c^{-1}q_a^{-1}q_0^{-1})
\Gamma_-'(yqq_c^{-1}q_a^{-1}q_0^{-1}) \\ \cdot 
\Gamma_-(yqq_a^{-1}q_0^{-1})
\Gamma_-'(yqq_0^{-1})
\end{split}
\]
so that 
\begin{equation}
\label{eqn:pyramid_normalized_product}
Z_{\text{pyramid}} = 
\left< \phi \left|
\cdots A'_+(q^2)A'_+(q) A'_+(1) A'_-(1) A'_-(q)A'_-(q^2) \cdots
\right| \phi \right>. 
\end{equation}
The commutation relation for these $A'$ operators is
\[
\begin{split}
A'_+(x)A'_-(y) = \frac{
(1+q_bxyq)(1+q_bq_cq_axyq)(1+q_b^{-1}xyq)(1+q_cxyq)
}{
(1-xyq)(1-q_bq_cxyq)(1-xyq)(1-q_cq_axyq)
}\\ \cdot 
\frac{
(1+q_c^{-1}xyq)(1+q_axyq)(1+(q_bq_cq_a)^{-1}xyq)(1-q_a^{-1}xyq)
}{
(1-(q_bq_c)^{-1}xyq)(1-xyq)(1-(q_cq_a)^{-1}xyq)(1-xyq)
} 
\\ \cdot A'_-(y)A'_+(x).
\end{split}
\]
Because of the mixed $\Gamma$ and $\Gamma'$ operators, some of the commutation factors now appear in the numerator.  We now move the $A'_+$ operators in (\ref{eqn:pyramid_normalized_product}) to the left of the expression, while moving the $A'_-$ operators to the right.   As in the proof of Theorem~\ref{thm:zn_result}, all of the $A'$ vanish, and we are left with the commutator
\[
Z_{\text{pyramid}} = 
\frac{
   M(1,q)^4 \widetilde{M}(q_bq_c, q)\widetilde{M}(q_aq_c, q)
}{
   \widetilde{M}(-q_a, q) \widetilde{M}(-q_b, q)
   \widetilde{M}(-q_c, q) \widetilde{M}(-q_aq_bq_c, q) 
},
\]
where $q = q_0q_aq_bq_c$.
\end{proof}
Note that this method gives an alternate proof of the result in~\cite{MYSELF}, when we specialize $q_0 = q_c$, $q_1 = q_a = q_b$.  

\section{Counting $\bbZ_2 \times \bbZ_2$--coloured 3D Young diagrams}

We will now prove Theorem~\ref{thm:z2z2_result}.  
Let us name the elements of $\bbZ_2 \times \bbZ_2$  $\{0, a, b, c\}$ as in the previous section, and recall the definition of $K_{\bbZ_2 \times \bbZ_2}$ from Theorem~\ref{thm:z2z2_result}.  Our set of indeterminates is $\mathcal{Q}=\{q_0, q_a, q_b, q_c\}$.  Let $q=q_0q_aq_bq_c$.  

\begin{figure}
\caption{Slicing a $\bbZ_2 \times \bbZ_2$--coloured 3D diagram}
\label{fig:z2z2_slice}
\begin{center}
\subfloat{
	\includegraphics[width=2in]{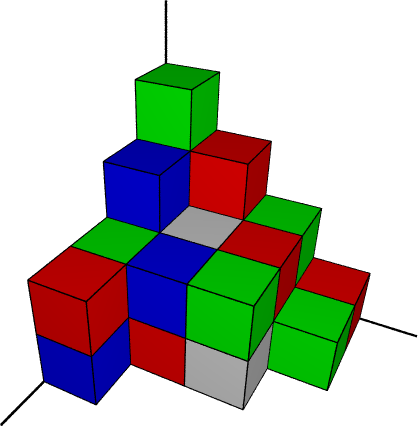}
}
\subfloat{
	\includegraphics[width=2in]{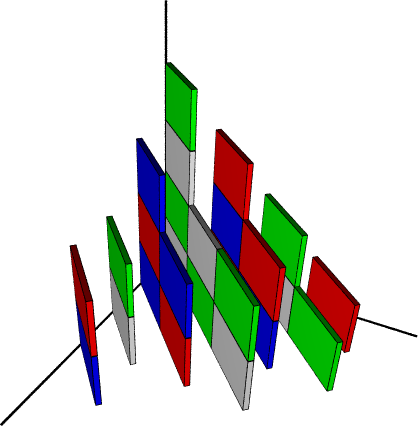}
}
\end{center}
\end{figure}

Before we proceed to compute this generating function, consider the $k$th diagonal slice $x-y = k$ of the positive octant.  Note that we have
\[
K_{\bbZ_2 \times \bbZ_2}(x, x+k, z) = (x-z)c + kb.
\]
In particular, the box $(x, x+k, z)$ is coloured $k \cdot b$ if $x \equiv z \pmod{2}$, and $k\cdot b + c$ otherwise.  In other words, each diagonal slice of $\pi$ is now coloured in a checkerboard fashion, whereas in the $\bbZ_n$ case, they were single--coloured (see Figure~\ref{fig:z2z2_slice}).  
We need to introduce a two--coloured weight function if we are to use vertex operators to compute $Z_{\bbZ_2 \times \bbZ_2}$.

\begin{definition} For $g,h \in \bbZ_2 \times \bbZ_2$, let
\[
Q_{gh} \lambda = 
q_g^{\#\{(i,j) \in \lambda \; | \; i \equiv j \pmod{2}\}} 
\cdot
q_h^{\#\{(i,j) \in \lambda \; | \; i \not \equiv j \pmod{2}\}} 
\cdot \lambda.
\]
\end{definition}
We may write down a vertex operator product which sweeps out a 3D diagram in diagonal slices, according to the $\bbZ_2 \times \bbZ_2$ colouring.  It is
\begin{equation}
\label{eqn:z2z2_product}
\begin{split}
Z_{\bbZ_2 \times \bbZ_2} = 
\left< \phi \right|
\cdots
Q_{ba}
\Gamma_+(1) Q_{0c} 
\Gamma_+(1) Q_{ba}
\Gamma_+(1) Q_{0c} \\ \cdot
\Gamma_-(1) Q_{ab}
\Gamma_-(1) Q_{0c}
\Gamma_-(1) Q_{ab}
\cdots
\left| \phi \right>
\end{split}
\end{equation}
Unfortunately, $\Gamma_{\pm}$ no longer commutes nicely with the $Q_{gh}$ operators, so our usual approach to computing with vertex operators fails here.  The problem is fundamental, and it does not appear that we can resolve it in a natural way.  We need a new idea.

However, there are two clues which tell us how to proceed.  The first clue is that the desired formula for $Z_{\bbZ_2 \times \bbZ_2}$ is very close to $Z_{\text{pyramid}}$, so it would suffice to prove the following:
\begin{lemma}
\label{lem:young_diagram_to_pyramid}
\[
Z_{\bbZ_2 \times \bbZ_2} = \widetilde{M}(q_aq_b, q) \cdot Z_{\text{pyramid}}.
\]
\end{lemma}
The second clue is that if we attempt to compute $\bbZ_{\text{pyramid}}$ by slicing along lines $x+y=k$, rather than $x-y=k$, then the slices of the pyramid partition are checkerboard--coloured as well!  See Figure~\ref{fig:pyramid_below}, which shows the colouring scheme from below.  The heavy black lines represent two edges of the pyramid, corresponding to prefixes of the words $(w_1v_1)^k$ and $(w_2v_2)^k$, so bricks which lie on these lines represent the corners of the slices.

So, using our checkerboard coloured weight operators, we can write down a different vertex operator product which still counts pyramid partitions.
\begin{lemma}
\label{lem:z2z2_pyramid_product}
\begin{equation*}
\begin{split}
Z_{\text{pyramid}} = 
\left< \phi \right|
\cdots
Q_{ba}
\Gamma_+'(1) Q_{0c}
\Gamma_+(1) Q_{ba}
\Gamma_+'(1) Q_{0c} \\ \cdot
\Gamma_-(1) Q_{ab}
\Gamma_-'(1) Q_{0c}
\Gamma_-(1) Q_{ab}
\cdots
\left| \phi \right>
\end{split}
\end{equation*}
\end{lemma}
\begin{proof}
One must check that the interlacing behaviour of the slices is correct, and that the correct weights are assigned to each slice.  This is similar to the proofs of Lemmas~\ref{lem:pyramid_slice} and~\ref{lem:pyramid_interlacing}.  
\end{proof}

Observe that~(\ref{eqn:z2z2_product}) is very similar to the product in Lemma~\ref{lem:z2z2_pyramid_product}, so we shall look for a way to transform $\Gamma_{\pm}(x)$ into $\Gamma_{\pm}'(x)$.

\begin{definition} Define
\[E_{\pm}(x) = \exp \sum_{k \geq 1} \frac{x^{2k}}{k} \alpha_{\pm 2k}.
\]
\end{definition}

\begin{lemma} The operators $E_{\pm}$ have the following properties:
\label{lem:e_commutation}
\begin{align*}
\Gamma_{\pm}(x) &= \Gamma'_{\pm}(x) E_{\pm}(x). \\
[E_{\pm}, \Gamma_{\pm}] &= 0. \\
E_{+}(x)\Gamma_-(y) &= \frac{1}{1-(xy)^2} \Gamma_-(y)E_{+}(x). \\
\Gamma_+(x)E_{-}(y) &= \frac{1}{1-(xy)^2} E_{-}(y)\Gamma_+(x). \\
\Gamma_+'(x)E_{-}(y) &= (1-(xy)^2) E_{-}(y)\Gamma_+'(x). 
\end{align*}
\end{lemma}

\begin{proof}
These are all simple applications of Corollaries~\ref{cor:CBH_I} and~\ref{cor:CBH_II} as well as (\ref{eqn:gamma_formula}) and~(\ref{eqn:second_gamma_formula}).
\end{proof}

In fact, unlike $\Gamma_{\pm}(x)$, $E_{\pm}(x)$ also commutes nicely with the checkerboard weight operators $Q_{gh}$.  

\begin{lemma}
\label{lem:e_commutes_with_checkerboard}
\begin{align*}
E_-(x) Q_{gh} &= Q_{gh} E_-(x\sqrt{q_gq_h});  \\
Q_{gh} E_+(x) &= E_+(x\sqrt{q_gq_h}) Q_{gh}.
\end{align*}
\end{lemma}

\begin{proof}
The operator $\alpha_{2n}$ has the effect of adding all possible border strips $R$ of length $2n$ to the boundary of a 2D Young diagram.  Since the length of the strips $R$ is even, any such $R$ has the same $Q_{gh}$--weight.  Indeed, 
\[Q_{gh} \cdot R = (q_gq_h)^n \cdot R = (\sqrt{q_gq_h})^{|R|} \cdot R,\]
It follows that
\[
\left(\sum_n \frac{x^{2n}}{n} \alpha_{-2n} \right) Q_{gh} \cdot \lambda
=
Q_{gh} \left(\sum_n \frac{(q\sqrt{q_gq_h})^{2n}}{n} \alpha_{-2n} \right) \cdot \lambda
\]
and thus
\[
E_-(x) Q_{gh} = Q_{gh} E_-(x\sqrt{q_gq_h}). 
\]
The case of $E_+$ is similar.
\end{proof}
Finally, we have the following property of $E_{\pm}$, inherited from the corresponding property of $\alpha_{\pm n}$:
\begin{align*}
\left< \phi \right| E_-(x) &= \left< \phi \right| \\
E_+(x) \left| \phi \right> &= \left| \phi \right>
\end{align*}

\emph{Proof of Lemma~\ref{lem:young_diagram_to_pyramid}.}
Let us alter the first line of~(\ref{eqn:z2z2_product}) by transforming half of the $\Gamma_+(1)$ operators into $\Gamma_+'(1)$ operators, 
\begin{align*}
&\!\!\!\!\!\!\!\!\!
\left<\phi\left|\prod_{i=1}^{\infty}\Gamma_+(1)Q_{ba} \Gamma_+(1) Q_{0c} \right|\right. \\
&= 
\left<\phi\left|\prod_{i=1}^{\infty}\Gamma_+(1)Q_{ba} \Gamma_+'(1) E_+(1)Q_{0c} \right|\right. \\
&=
\left<\phi\left|\prod_{i=1}^{\infty}\Gamma_+(1)Q_{ba} \Gamma_+'(1)Q_{0c} \right|\right. \cdot
\left|\prod_{i=0}^{\infty}E_+(q^i)E_+(Q^i\sqrt{q_bq_a})\right|.
\end{align*}
Now, continue to move the $E_+$ term to the right through the second line of~(\ref{eqn:z2z2_product}).  We have
\begin{align*}
&\!\!\!\!\!\!\!\!\!
\left|\prod_{i=0}^{\infty}E_+(q^i)E_+(Q^i\sqrt{q_bq_a})\right|\cdot
\left|\left.
\prod_{i=1}^{\infty}\Gamma_-(1)Q_{ab} \Gamma_-(1) Q_{0c}
\right| \phi \right> \\
&= C \cdot
\left|\left.
\prod_{i=1}^{\infty}\Gamma_-(1)Q_{ab} \Gamma_-(1) Q_{0c}
\right|  \cdot
\left|\prod_{i=0}^{\infty}E_+(q^i)E_+(Q^i\sqrt{q_bq_a})\right|
\phi \right> \\
&= C \cdot
\left|\left.
\prod_{i=1}^{\infty}\Gamma_-(1)Q_{ab} \Gamma_-(1) Q_{0c}
\right|  
\phi \right>, \\
\end{align*}
where $C = M(1,q) M(q_b^{-1}q_c^{-1}, q)$ is the product of the commutators generated by Lemma~\ref{lem:e_commutation}.  Next, change half of the $\Gamma_-$ to $\Gamma_-'$ in the above expression,
\[
\left|\left.
\prod_{i=1}^{\infty}\Gamma_-(1)Q_{ab} \Gamma_-'(1)E_-(1) Q_{0c}
\right|  
\phi \right> 
\]
and commute them out to the left.  This time, we pick up the multiplicative factor 
\[
\frac{M(q_aq_b, q)}{M(1,q)},
\]
and the factors $M(1,q)$ cancel.  We have shown that 
\[
Z_{\bbZ_2 \times \bbZ_2} = Z_{\text{pyramid}} \cdot \widetilde{M}(q_aq_b).
\]
so Lemma~\ref{lem:young_diagram_to_pyramid} and Theorem~\ref{thm:z2z2_result} are now proven.$\hfill \square$

\section{Future work}

It would obviously be wonderful to have a combinatorial proof of these identities; such a proof might be analagous to the ``$n$--quotient'' on 2D Young diagrams, which decomposes a Young diagram into $n$ smaller Young diagrams and an $n$--core.   The authors suspect, however, that such a proof would be rather difficult to find.   One indication of this is that there are formulae for 3D Young tableaux which fit inside an $A \times B \times C$ box, but computational evidence suggests that there is no such nice formula for $\bbZ_2 \times \bbZ_2$--coloured partitions.

One could attempt to compute the Don\-ald\-son--Thom\-as partition functions of arbitrary toric Calabi-Yau orbifolds.  To this aim, it should be possible to develop an orbifold version of the topological vertex formalism following~\cite{MNOP}; this is a work in progress.  
It would also be interesting to try to extend Szendr\H{o}i's work~\cite{NC_DT} in noncommutative Don\-ald\-son--Thom\-as theory using the results of this paper.  

One box counting problem which is of great interest is to take $G=\bbZ_3$ and the colouring $K$ given by
\[
K(1,0,0) = K(0,1,0) = K(0,0,1) = 1.
\]
However, this problem appears to be rather difficult.  The group representation does not naturally embed into $SO(3)$ or $SU(2)$, so Don\-ald\-son--Thom\-as theory does not generate any conjectures as to what the answer might look like.  Indeed, Kenyon~\cite{KENYON} conjectures that there is no nice product formula in this example.

One unifying theme between 3D diagrams and pyramid partitions is quivers: both objects arise from a quiver path algebra modulo an ideal generated by a superpotential~\cite{NC_DT}.  Perhaps one can extend the methods to other quivers and superpotentials.

However, the most intriguing direction for future work is simply to try to understand these proofs more fully.  The reader may perhaps have noticed that the appearance of pyramid partitions seems somewhat unmotivated.  Undoubtedly, there is some underlying geometric or rep\-re\-sen\-ta\-tion-theoretic reason why these product formulae exist.

\appendix
\section{Don\-ald\-son-Thom\-as theory of $\cnums^{3}/G$ and its
crepant resolution (by Dr. Jim Bryan)}
\label{app: DT theory of C3/G}

\subsection{Review of Don\-ald\-son-Thom\-as theory}~

Don\-ald\-son-Thom\-as theory, in its incarnation due Mau\-lik, Okoun\-kov,
Nek\-rasov, and Pand\-hari\-pande, constructs subtle integer valued
deformation invariants of a threefold $X$ out of the Hilbert scheme of
subschemes of $X$. If $X$ is a Calabi-Yau threefold, i.e., $K_{X} $ is
trivial, then this invariant has a simple formulation due to
Behrend. It is given by the weighted topological Euler characteristic
of the Hilbert scheme where the weighting is by $\nu $, an integer
valued constructible function which is canonically associated to any
scheme \cite{Behrend}.

Let $X$ be a (not necessarily compact) threefold with trivial
canonical class. Let $I_{n} (X,\beta )$ be the Hilbert scheme of
subschemes $Z\subset X$ having proper support of dimension less than
or equal to one and with $[Z]=\beta \in H_{2} (X)$ and $n=\chi (\O
_{Z})$. We define the Don\-ald\-son-Thom\-as invariant $N_{\beta }^{n} (X)$
to be
\begin{align*}
N_{\beta }^{n} (X) 
&=e (I_{n} (X,\beta ), \nu )\\
&= \sum_{k\in \znums} k e\left( \nu^{-1}(k)\right) \\
\end{align*}
where $e (\cdot )$ denotes topological Euler
characteristic and $\nu $ is Behrend's constructible function .

The invariants are assembled into the partition function $Z^{DT}_{X}$
as follows. Let $C_{1},\dotsc ,C_{l}$ be a basis for $H_{2} (X,\znums
)$ such that any effective curve class $\beta $ is given by
$d_{1}C_{1}+\dotsb +d_{l}C_{l}$ with $d_{i}\geq 0$. Let $v_{1},\dotsc
,v_{l}$ be corresponding variables and let $v^{\beta
}=v_{1}^{d_{1}}\dotsb v_{l}^{d_{l}}$. The Don\-ald\-son-Thom\-as partition
function of $X$ is defined by
\[
Z^{DT}_{X} (v,q) = \sum _{\beta \in H_{2} (X,\znums )} \sum _{n\in
\znums } N_{\beta }^{n} (X)v^{\beta } q^{n}.
\]
Define the reduced partition function by
\begin{align*}
Z^{DT}_{X} (v,q)' &= \frac{Z_{X}^{DT} (v,q)}{Z^{DT}_{X} (0,q)}\\
&={M (1,q)^{-e (X)}}Z^{DT}_{X} (v,q)
\end{align*}
where the second equality is a theorem proved by
\cite{Behrend-Fantechi,Levine-Pandharipande, Li}.

Maulik, Nekrasov, Okounkov, and Pandharipande conjecture that
Don\-ald\-son-Thom\-as theory is equivalent to Gro\-mov-Wit\-ten theory.  We
assemble $GW^{g}_{\beta } (X)$, the genus $g$ Gro\-mov-Wit\-ten invariants
of non-zero degree $\beta $ into the reduced Gro\-mov-Wit\-ten partition
function as follows:
\[
Z^{GW}_{X} (v,\lambda )' =\exp\left(\sum _{\beta \neq 0}\sum
_{g=0}^{\infty } GW_{\beta }^{g} (X)v^{\beta }\lambda ^{2g-2} \right).
\]

\begin{conjecture}\cite{MNOP1}
Under the change of variables $q=-e^{i\lambda }$ the reduced partition
functions of Don\-ald\-son-Thom\-as and Gro\-mov-Wit\-ten theory are equal:
\[
Z^{DT}_{X} (v,q)' = Z^{GW}_{X} (v,\lambda )'.
\]
\end{conjecture}

This conjecture has been proven in the case where $X$ is a toric
local surface \cite{MNOP1} and when $X$ is a local curve
\cite{Bryan-Pandharipande-localcurves,Okounkov-Pandharipande-dtlc}. 

\smallskip

\subsection{Orbifold Don\-ald\-son-Thom\-as theory of $[\cnums ^{3}/G]$}~

Extending Do\-nald\-son-Tho\-mas theory to the case of three dimensional
orbifolds is expected to be routine since the Hilbert scheme of
substacks of a Deligne-Mumford stack has been constructed by Olsson
and Starr \cite{Olsson-Starr}, although it isn't clear how best to
choose the discrete data in general.

The orbifolds that we consider are simple enough that we can identify
the Hilbert scheme directly. Let $G$ be a finite subgroup of $SU
(3)$. A substack of $[\cnums ^{3}/G]$ may be regarded as a
$G$-invariant subscheme of $\cnums ^{3}$, and consequently we can
regard the Hilbert scheme of $[\cnums ^{3}/G]$ as a subset of the
Hilbert scheme of $\cnums ^{3}$. Since we require our substacks to
have proper support, we need only consider zero dimensional subschemes
of $\cnums ^{3}$. For any $G$-rep\-re\-sen\-ta\-tion $R$ of dimension $d$ we
identify $\Hilb ^{R} ([\cnums ^{3}/G])\subset \Hilb ^{d} (\cnums
^{3})$ as follows:
\[
\Hilb ^{R} ([\cnums ^{3}/G]) = \left\{Z\subset \cnums ^{3}:
\text{ $Z$ is $G$-invariant with
$H^{0} (\O _{Z})=R$} \right\}.
\]

This Hilbert scheme has a symmetric perfect obstruction theory induced
by the $G$-invariant part of the of the perfect obstruction theory on
$\Hilb ^{d} (\cnums ^{3})=I_{d} (\cnums ^{3},0)$. However, we do not
need this construction since we can define the Don\-ald\-son-Thom\-as
invariants directly using Behrend's constructible function.

\begin{definition}
The Don\-ald\-son-Thom\-as invariants of $[\cnums ^{3}/G]$ are indexed by
representations of $G$ and are given by the Euler characteristics of
the Hilbert schemes, weighted by Behrend's $\nu  $ function:
\[
N^{R} (\cnums ^{3}/G) = e (\Hilb ^{R} ([\cnums ^{3}/G]),\nu ).
\]
\end{definition}

Let $q_{0},\dotsc ,q_{r}$ be variables corresponding to $R_{0},\dotsc
,R_{r}$, the irreducible representations of $G$. For a representation
$R=d_{0}R_{0}+\dotsb +d_{r}R_{r}$, let $q^{R} $ denote
$q_{0}^{d_{0}}\dotsb q_{r}^{d_{r}}$. We define the orbifold
Don\-ald\-son-Thom\-as partition function by
\[
Z^{DT}_{\cnums ^{3}/G} (q_{0},\dotsc ,q_{l}) = \sum _{R } N^{R}
(\cnums ^{3}/G) q^{R}
\]
where $R$ runs over all representations of $G$.

We now restrict our attention to groups $G$ which are subgroups of $SO
(3)\subset SU (3)$ and are Abelian. Finite subgroups of $SO (3)$ admit
an ADE classification. They are the cyclic groups, the dihedral
groups, and the platonic groups. The only Abelian groups from this
list are the cyclic groups $\znums _{n}$ and the Klein 4-group
$\ZtwoZtwo $. The action of $k\in \znums _{n}$ on $\cnums ^{3}$ is
given by
\[
k (x,y,z) =(\omega  ^{k}x,\omega  ^{-k}y,z)
\]
where $\omega =\exp\left(\frac{2\pi i}{n} \right)$. The action of
$\ZtwoZtwo =\{0,a,b,c \}$ on $\cnums ^{3}$ is given by
\begin{align*}
a (x,y,z)&= (x,-y,-z),\\
b (x,y,z)&= (-x, y,-z),\\
c (x,y,z)&= (-x,-y,z).
\end{align*}

As in the introduction, we choose an isomorphism $\psi $ of the group
of representations $\hat{G}$ with $G$. Explicitly, we identify $1\in
\znums _{n}$ with $L$, the representation of $\znums _{n}$ where $1\in
\znums _{n}$ acts by multiplication by $\exp\left(\frac{2\pi i}{n}
\right)$. For $\ZtwoZtwo =\{0,a,b,c \}$ we identify $a,$ $b$, and $c$
with the representations $\alpha $, $\beta $, and $\gamma $ given by
the action on the $x$, $y$, and $z$ coordinates of $\cnums ^{3}$
respectively.

\begin{theorem}\label{thm: change of signs on variables}
Let $q_{k}$ be the variable corresponding to the group element $k\in
\znums _{n}$ and the character $L^{k}$. Then
\[
Z^{DT}_{\cnums ^{3}/\znums _{n}} (q_{0},\dotsc ,q_{n-1}) = Z_{\znums
_{n}} (-q_{0},q_{1},\dotsc ,q_{n-1})
\]
where $Z_{\znums _{n}}$ is the $\znums _{n}$-coloured 3D diagram
partition function introduced and computed in the main body of the
paper (Theorem~\ref{thm:zn_result}). 

Let $\{q_{0},q_{a},q_{b},q_{c} \}$ be variables corresponding to
$\{0,a,b,c \}$, the group elements of $\ZtwoZtwo $, and $\{1,\alpha
,\beta ,\gamma \}$, the characters of $\ZtwoZtwo $. Then
\[
Z_{\cnums ^{3}/\ZtwoZtwo }^{DT} (q_{0},q_{a},q_{b},q_{c}) =
Z_{\ZtwoZtwo } (q_{0},-q_{a},-q_{b},-q_{c})
\]
where $Z_{\ZtwoZtwo }$ is the $\ZtwoZtwo $-coloured 3D diagram
partition function introduced and computed in the main body of the
paper (Theorem~\ref{thm:z2z2_result}). 
\end{theorem}

\textsc{Proof:} Let $G$ be $\znums _{n}$ or $\ZtwoZtwo $ and let
$T\subset \left(\cnums ^{\times } \right)^{3}$ be the subtorus with
$t_{1}t_{2}t_{3}=1$. The action of $T$ on $\cnums ^{3}$ commutes with
the action of $G$ and hence defines a $T$-action on $[\cnums ^{3}/G] $
and on $\Hilb ^{R} (\cnums ^{3}/G)$. The fixed points of $T$ in $\Hilb
^{R} (\cnums /G)\subset \Hilb ^{\dim R} (\cnums ^{3})$ are isolated,
even infinitesimally \cite[Lemma~4.1]{Behrend-Fantechi}, and they
correspond to monomial ideals in $\cnums [x,y,z]$. The monomial ideals
in turn correspond to 3D Young diagrams $\pi $ where if $Z$ denotes
the $T$-fixed subscheme of $\cnums ^{3}$, then
\[
H^{0} (\O _{Z}) = \sum _{(i,j,k)\in \pi } t_{1}^{i}t_{2}^{j}t_{3}^{k}
\]
as a $T$-rep\-re\-sen\-ta\-tion viewed as a polynomial in $t_{1},t_{2},t_{3}$
modulo the relation $t_{1}t_{2}t_{3}=1$. Following \cite{MNOP1}, we
adopt the notation
\[
Q_{\pi } = \sum _{(i,j,k)\in \pi } t_{1}^{i}t_{2}^{j}t_{3}^{k}.
\]
By \cite[Prop.~3.3]{Behrend-Fantechi}, the $\nu $-weighted Euler
characteristic of $\Hilb ^{R} (\cnums ^{3}/G)$ is given by a sum over
the $T$-fixed points, counted with sign given by the parity of the
dimension of the Zariski tangent space of $\Hilb ^{R} (\cnums ^{3}/G)$
at a fixed point corresponding to 3D diagram $\pi $. Hence both the
Don\-ald\-son-Thom\-as and the diagram partition functions are given by a
sum over 3D diagrams, weighted, up to a sign, by the same
variables. Thus our main task is to determine the sign.
 
Let $\pi $ be a 3D diagram having $N=|\pi |$ boxes and having $|\pi
|_{g}$ boxes of colour $g\in G$. Let $T_{\pi }$ denote the Zariski
tangent space of $\Hilb ^{N} (\cnums ^{3})$ at the subscheme
corresponding to $\pi $. Let
\[
\left(T_{\pi } \right)^{0}\subset  T_{\pi }
\]
be the Zariski tangent space of $\Hilb ^{R} ([\cnums ^{3}/G])\subset
\Hilb ^{N} (\cnums ^{3})$ at the same point. $T_{\pi }$ can be
regarded as both a $T$-rep\-re\-sen\-ta\-tion and a $G$-rep\-re\-sen\-ta\-tion.
$\left(T_{\pi } \right)^{0}$ is given by the $G$-invariant subspace of
$T_{\pi }$.

The difference of $T_{\pi }$ and its dual $T_{\pi }^{\vee }$, regarded
as a virtual $\left(\cnums ^{\times } \right)^{3}$-rep\-re\-sen\-ta\-tion, is
computed in \cite[equation~(13)]{MNOP1} and given by
\[
T_{\pi }-T_{\pi }^{\vee } = Q_{\pi } - \frac{\overline{Q}_{\pi
}}{t_{1}t_{2}t_{3}} + Q_{\pi }\overline{Q}_{\pi } \frac{(1-t_{1})
(1-t_{2}) (1-t_{3})}{t_{1}t_{2}t_{3}}
\]
where 
\[
\overline{Q}_{\pi } (t_{1},t_{2},t_{3})=Q_{\pi }
(t_{1}^{-1},t_{2}^{-1},t_{3}^{-1}).
\]
Using the relation
$t_{1}t_{2}t_{3}=1$ to eliminate $t_{3}$ from the above expression, we
can regard $T_{\pi } -T^{\vee }_{\pi }$ as an element in 
\[
R (T)\cong \znums [t_{1},t_{2},t_{1}^{-1},t_{2}^{-1}],
\]
the virtual representation ring of $T$.

Following \cite{MNOP1}, we let 
\[
V_{\pi } = Q_{\pi } + Q_{\pi }\overline{Q}_{\pi } (1-t_{1})
(1-t_{2})t_{1}^{-1}t_{2}^{-1}
\]
which satisfies the easily verified equation 
\begin{equation}\label{eqn: T-Tv=V-Vv}
T_{\pi }-T^{\vee }_{\pi } = V_{\pi }- V^{\vee }_{\pi }
\end{equation}
in $R (T)$, and also has the crucial property that the constant term of
$V_{\pi }$ is even \cite[Lemma~10]{MNOP1}. These facts allow us to use
$V_{\pi }$ as a surrogate for $T_{\pi }$ when computing the parity of
the dimension:

\begin{lemma}
Let $\left(T_{\pi } \right)^{0}$ and $\left(V_{\pi } \right)^{0}$
denote the $G$-invariant part of $T_{\pi }$ and $V_{\pi }$
respectively, then
\[
\dim \left(T_{\pi } \right)^{0}  \equiv \vdim\left(V_{\pi } \right)^{0} \mod 2.
\]
\end{lemma}
\textsc{Proof:} From equation~\eqref{eqn: T-Tv=V-Vv} we see that
$T_{\pi }-V_{\pi }$ is self-dual. Thus all non-constant monomials
occur in pairs of the form $a_{ij}
(t_{1}^{i}t_{2}^{j}+t_{1}^{-i}t_{2}^{-j})$. Moreover, the constant
term of $V_{\pi }$ is even \cite[Lemma~10]{MNOP1} and the constant
term of $T_{\pi }$ is zero \cite[Lemma~4.1]{Behrend-Fantechi}. Thus we
have
\[
\vdim \left(T_{\pi }-V_{\pi } \right) \equiv 0 \mod 2.
\]
Indeed, the above argument shows that if we restrict $T_{\pi }-V_{\pi
}$ to any self-dual collection of weights, the virtual dimension will
be even. In particular, the $G$-invariant part of $T_{\pi }-V_{\pi }$
has even virtual dimension, which proves the lemma. \qed

To compute the parity of the $G$-invariant part of $V_{\pi }$, we work
in the representation ring of $G$ with mod 2 coefficients. The restriction map
\[
R (T)\cong  \znums [t_{1},t_{2},t_{1}^{-1},t_{2}^{-1}] \to R_{\znums _{2}} (G)
\]
is explicitly given by
\[
(t_{1},t_{2})\mapsto (L,L^{-1})
\]
in the case where $G=\znums _{n}$, and by 
\[
(t_{1},t_{2})\mapsto (\alpha ,\beta )
\]
in the case where $G=\ZtwoZtwo $. 

For any $W\in R_{\znums _{2}} (G)$ and any irreducible representation
$\zeta $, let $[W]_{\zeta }\in \znums _{2}$ denote the coefficient of
$\zeta $ in $W$. We compute $[V_{\pi }]_{1}$ in
\[
R_{\znums _{2}} (\znums _{n})=\znums _{2}[L ]/ (L ^{n}-1)
\]
as follows.
\begin{align*}
[V_{\pi }]_{1} &=\left[Q_{\pi }+Q_{\pi }\overline{Q}_{\pi } (1-L ) (1-L ^{-1}) \right]_{1}\\
&=[Q_{\pi }]_{1}+\left[Q_{\pi }\overline{Q}_{\pi } (L +L ^{-1}) \right]_{1}\\
&=[Q_{\pi }]_{1} + \left[Q_{\pi }\overline{Q}_{\pi } \right]_{L ^{-1}} +\left[Q_{\pi }\overline{Q}_{\pi } \right]_{L }\\
&=[Q_{\pi }]_{1}\\
&=|\pi |_{0} \mod 2
\end{align*}
Since $[V_{\pi }]_{1}$ is equal to the dimension of the $\znums
_{n}$-invariant part of $T_{\pi }$ modulo 2, the 3D diagram $\pi $ is
counted with sign $(-1)^{|\pi |_{0}}$ in the Don\-ald\-son-Thom\-as
partition function of $\cnums ^{3}/\znums _{n}$. This proves first
part of Theorem~\ref{thm: change of signs on variables}.

We now compute $[V_{\pi }]_{1}$ in 
\[
R_{\znums _{2}} (\ZtwoZtwo ) = \znums _{2}[\alpha ,\beta ]/ (\alpha
^{2}-1,\beta ^{2}-1).
\]
We use the fact that in this ring, the square of an arbitrary element
is equal to the sum of its coefficients:
\[
(n_{1}+n_{2}\alpha +n_{3}\beta +n_{4}\alpha \beta )^{2} = n_{1}+n_{2}+n_{3}+n_{4} ,
\]
and we compute as follows.
\begin{align*}
[V_{\pi }]_{1} &=\left[Q_{\pi }+Q_{\pi }\overline{Q}_{\pi } (1-\alpha  ) (1-\beta )\alpha \beta  \right]_{1}\\
&=[Q_{\pi }]_{1}+\left[Q_{\pi }^{2}(1+\alpha +\beta +\alpha \beta ) \right]_{1}\\
&=[Q_{\pi }]_{1} + \left[|\pi | (1+\alpha +\beta +\alpha \beta ) \right]_{1}\\
&=|\pi |_{0} +|\pi | \mod 2\\
&=|\pi |_{a}+|\pi |_{b}+|\pi |_{c} \mod 2.
\end{align*}
Since $[V_{\pi }]_{1}$ is equal to the dimension of the $\ZtwoZtwo
$-invariant part of $T_{\pi }$ modulo 2, the 3D diagram $\pi $ is
counted with sign $(-1)^{|\pi |_{a}+|\pi |_{b}+|\pi |_{c}}$ in the
Don\-ald\-son-Thom\-as partition function of $\cnums ^{3}/\ZtwoZtwo $. This
proves the remaining part of Theorem~\ref{thm: change of signs on
variables} and so the proof of Theorem is complete.\qed

\begin{remark}
For any finite Abelian subgroup $G\subset SU (3)$, the
Don\-ald\-son-Thom\-as invariants of $\cnums ^{3}/G$ are given by a signed
count of boxes coloured by $G$. However, it is not always true that
this sign is obtained by simply changing the signs of some of the
variables. For example, consider the case of $G=\znums _{3}$ acting on
$\cnums ^{3}$ with equal weights on all three factors. The sign
associated to a 3D partition $\pi $ can be computed by the methods of
this appendix and is given by $(-1)^{\sigma}$, where
\[
\sigma = 
|\pi|_1 + |\pi|_2 + |\pi|_0 |\pi|_1 + |\pi|_0 |\pi|_2 + |\pi|_1 |\pi|_2.
\]
Thus the coloured 3D diagram partition function and the
Don\-ald\-son-Thom\-as partition function are not related in an obvious way.
\end{remark}

\smallskip

\subsection{The Don\-ald\-son-Thom\-as crepant resolution conjecture}~

A well known principle in physics asserts that string theory on a
Calabi-Yau orbifold $X$ is equivalent to string theory on any
crepant resolution $Y\to X$. Consequently, it is expected that
mathematical counterparts of string theory, such as Gro\-mov-Wit\-ten
theory or Don\-ald\-son-Thom\-as theory, should be equivalent on $X$ and
$Y$. Precise formulations of these equivalences are known as crepant
resolution conjectures. The crepant resolution conjecture in
Gro\-mov-Wit\-ten theory goes back to Ruan, and has recently undergone
successive refinements \cite{Ruan-crepant, Bryan-Graber, CCIT,
Coates-Ruan}.

In this section we formulate a crepant resolution conjecture for
Don\-ald\-son-Thom\-as theory. Our conjecture has somewhat limited scope: we
stick to the ``local case'' where $X$ is of the form $[\cnums
^{3}/G]$, and (for reasons explained below) we impose the \emph{hard
Lefschetz condition} \cite[Defn~1.1]{Bryan-Graber}, which implies
\cite{Bryan-Gholampour3} that $G$ is a finite subgroup of either $SU
(2)\subset SU (3)$ or $SO (3)\subset SU (3)$.

The most straightforward formulation of the crepant resolution
conjecture in Don\-ald\-son-Thom\-as theory posits that the partition
functions of the orbifold and its resolution are equal after some
natural change of variables. For the orbifold $[\cnums ^{3}/G]$, we
saw in the previous section that the partition function has variables
naturally indexed by irreducible $G$-rep\-re\-sen\-ta\-tions. By the
classical McKay correspondence, the crepant resolution $Y_{G}\to
\cnums ^{3}/G$ given by the $G$-Hilbert scheme has a basis of $H_{*}
(Y_{G})$ also labelled by irreducible $G$-rep\-re\-sen\-ta\-tions
\cite{BKR,McKay}. However, the variables of the Don\-ald\-son-Thom\-as
partition function of $Y_{G}$ correspond to a basis of $H_{0}
(Y_{G})\oplus H_{2} (Y_{G})$. So in order to get the number of
variables of $Z^{DT}_{Y_{G}}$ and $Z^{DT}_{\cnums ^{3}/G}$ to match,
we need
\[
H_{*} (Y_{G}) = H_{0} (Y_{G})\oplus H_{2} (Y_{G}).
\]
This occurs if and only if $Y_{G}\to \cnums ^{3}/G$ is a semi-small
resolution. This condition is equivalent to the orbifold satisfying
the hard Lefschetz condition.

\begin{conjecture}\label{conj: local DT CRC}
Let $X$ be a local, 3 dimensional, Calabi-Yau orbifold satisfying the
hard Lefschetz condition, namely, $X=X_{G}=[\cnums ^{3}/G]$ where $G$
is a finite subgroup of either $SU (2)\subset SU (3)$ or $SO
(3)\subset SU (3)$. 

Let $q_{0},q_{1},\dotsc ,q_{l}$ be variables
corresponding to the irreducible $G$-rep\-re\-sen\-ta\-tions
$R_{0},R_{1},\dotsc ,R_{l}$ where $R_{0}$ is the trivial
representation. Let $Y_{G}\to X_{G}$ be the crepant resolution given
by the $G$-Hilbert scheme and let $v_{1},\dotsc ,v_{l}$ be the
variables corresponding to the basis of curve classes in $Y_{G}$
labelled by the non-trivial $G$-rep\-re\-sen\-ta\-tions $R_{1},\dotsc ,R_{l}$.

Then the Don\-ald\-son-Thom\-as partition functions of $Y_{G}$ and $X_{G}$
are related by the formula
\[
Z^{DT}_{X_{G}} (q_{0},\dotsc ,q_{l}) = M (1,q)^{-e (Y_{G})}
Z^{DT}_{Y_{G}} (q,v_{1},\dotsc ,v_{l}) Z^{DT}_{Y_{G}}
(q,v_{1}^{-1},\dotsc ,v_{l}^{-1})
\]
under the identification of the variables
\begin{align*}
v_{i}&=q_{i} \quad \text{ for }\quad  i=1,\dotsc ,l,\\
q&=q^{R_{reg}}\\
&=q_{0}^{\dim R_{0}}\dotsb q_{l}^{\dim R_{l}}.
\end{align*}
\end{conjecture}

\begin{proposition}\label{prop: DT CRC holds for Zn and Z2xZ2}
Conjecture~\ref{conj: local DT CRC} holds for $G$ Abelian, namely for
$G=\znums _{n}$ or $G=\ZtwoZtwo $.
\end{proposition}

\begin{remark}
Szendr\H{o}i proved \cite{NC_DT} that a similar relationship holds
between the Don\-ald\-son-Thom\-as partition functions of the non-commutative
conifold singularity and its small resolution.
\end{remark}

\begin{remark}
The Gro\-mov-Wit\-ten partition function of $Y_{G}$ has been computed for
all $G$ in $SU (2)$ or $SO (3)$ in \cite{Bryan-Gholampour3} (see also
Remark~\ref{rem: choice of weights in equiv GW theory}). This
provides, via the MNOP conjecture, a prediction for $Z^{DT}_{Y_{G}}$
and hence our conjecture~\ref{conj: local DT CRC} gives a prediction
for $Z^{DT}_{\cnums ^{3}/G}$ which can be tested term by
term. Verification of this prediction for terms of low order has been
obtained by D.~Steinberg in the case where $G$ is the quaternion 8
group.
\end{remark}

In light of Theorems~\ref{thm:zn_result}, \ref{thm:z2z2_result}, and
\ref{thm: change of signs on variables}, Proposition~\ref{prop: DT CRC
holds for Zn and Z2xZ2} is equivalent to
Theorem~\ref{thm:resolution_dt_partition_functions} which we prove
here.

\subsubsection{Proof of Proposition~\ref{prop: DT CRC holds for Zn and
Z2xZ2} / Theorem~\ref{thm:resolution_dt_partition_functions}:}

Since $G$ is Abelian, $Y_{G}$ is toric and so via \cite[Theorems 2 and
3]{MNOP1}, the reduced Don\-ald\-son-Thom\-as partition function of $Y_{G}$
is equal to the reduced Gro\-mov-Wit\-ten partition function of $Y_{G}$
after the change of variables $q=-e^{i\lambda }$. Thus it suffices to
compute the Gro\-mov-Wit\-ten partition function of $Y_{G}$\footnote{In
\cite{MNOP1} it is shown that the reduced Don\-ald\-son-Thom\-as partition
function of a toric Calabi-Yau threefold can be computed via the
topological vertex formalism. In general, the topological vertex
formalism has been proven to compute the Gro\-mov-Wit\-ten partition
function only in the ``two-leg'' case. While $Y_{\znums _{n}}$ is a
local surface and can be computed with two-leg vertices, $Y_{\ZtwoZtwo
}$ has the geometry of the closed topological vertex \cite{Bryan-Karp}
and requires a three-leg vertex. However, in this case, the invariants
have been computed by both the vertex formalism as well as by
localization and have been shown to agree \cite{Ka-Liu-Ma}. Thus we
know that the Gro\-mov-Wit\-ten/Don\-ald\-son-Thom\-as correspondence holds for
both $Y_{\znums _{n}}$ and $Y_{\ZtwoZtwo }$.}. The pithiest way to
encode the Gro\-mov-Wit\-ten invariants is in terms of Gopakumar-Vafa
invariants, or so called BPS state counts. It is well know that each
genus zero BPS state count $n_{\beta }^{0}$ contributes a factor of $M
(v^{\beta },-e^{i\lambda })^{-n_{\beta }^{0}}$ to the Gro\-mov-Wit\-ten
partition function (see for example the proof of Theorem~3.1 in
\cite{Behrend-Bryan}). Thus the content of
Theorem~\ref{thm:resolution_dt_partition_functions} is that $Y_{\znums
_{n}}$ has genus 0 Gopakumar-Vafa invariants occurring in the classes
$C_{a}+\dotsb +C_{b}$ for $0<a\leq b<n$ with value -1, and that
$Y_{\ZtwoZtwo }$ has genus 0 Gopakumar-Vafa invariants occurring in
the classes $C_{a}$, $C_{b}$, $C_{c}$, and $C_{a}+C_{b}+C_{c}$ with
value 1 and in the classes $C_{a}+C_{b}$, $C_{a}+C_{c}$, and
$C_{b}+C_{c}$ with value -1. Moreover, all other Gopakumar-Vafa
invariants are zero. These assertions are proved in \cite{Ka-Liu-Ma}:
the case of $Y_{\ZtwoZtwo }$ is Corollary~16 and Proposition~19 and
the case of $Y_{\znums _{n}}$ is Proposition~12. \qed

\begin{remark}\label{rem: choice of weights in equiv GW theory}
The Gro\-mov-Wit\-ten and Don\-ald\-son-Thom\-as theories of $Y_{G}$ are
equivariant theories and so in general depend on the choice of the
torus action. In this paper, we have assumed that the torus is chosen
to act trivially on the canonical class. This choice is required to
apply the topological vertex formalism as we have done in the above
proof. We warn the reader that the computation of the Gro\-mov-Wit\-ten
invariants of $Y_{G}$ for general $G\subset SO (3)$ done in
\cite{Bryan-Gholampour3} is done using the $\cnums ^{\times }$ action
induced from the \emph{diagonal} action on $\cnums ^{3}/G$. This does
not change which classes carry Gopakumar-Vafa invariants, but it can
change the values of the invariants in those curve classes that admit
deformations to infinity.
\end{remark}

\bibliography{orbifold_dt}

\end{document}

%% file: add_ribbon.pdf_t
\begin{picture}(0,0)%
\includegraphics{add_ribbon.pdf}%
\end{picture}%
\setlength{\unitlength}{3947sp}%
\begingroup\makeatletter\ifx\SetFigFont\undefined%
\gdef\SetFigFont#1#2#3#4#5{%
  \reset@font\fontsize{#1}{#2pt}%
  \fontfamily{#3}\fontseries{#4}\fontshape{#5}%
  \selectfont}%
\fi\endgroup%
\begin{picture}(4824,2574)(1039,1277)
\put(2101,2414){\makebox(0,0)[lb]{\smash{{\SetFigFont{12}{14.4}{\rmdefault}{\mddefault}{\updefault}{\color[rgb]{0,0,0}$h=3$}%
}}}}
\put(3901,2264){\makebox(0,0)[lb]{\smash{{\SetFigFont{12}{14.4}{\rmdefault}{\mddefault}{\updefault}{\color[rgb]{0,0,0}$h=2$}%
}}}}
\put(4201,1364){\makebox(0,0)[lb]{\smash{{\SetFigFont{12}{14.4}{\rmdefault}{\mddefault}{\updefault}{\color[rgb]{0,0,0}$+$}%
}}}}
\put(2926,1364){\makebox(0,0)[lb]{\smash{{\SetFigFont{12}{14.4}{\rmdefault}{\mddefault}{\updefault}{\color[rgb]{0,0,0}$-$}%
}}}}
\put(1501,1364){\makebox(0,0)[lb]{\smash{{\SetFigFont{12}{14.4}{\rmdefault}{\mddefault}{\updefault}{\color[rgb]{0,0,0}$+$}%
}}}}
\put(5326,2639){\makebox(0,0)[lb]{\smash{{\SetFigFont{12}{14.4}{\rmdefault}{\mddefault}{\updefault}{\color[rgb]{0,0,0}$h=1$}%
}}}}
\put(3526,2864){\makebox(0,0)[lb]{\smash{{\SetFigFont{12}{14.4}{\rmdefault}{\mddefault}{\updefault}{\color[rgb]{0,0,0}$\alpha_{-3}$}%
}}}}
\end{picture}%

%% file: quiver.pdf_t
\begin{picture}(0,0)%
\includegraphics{quiver.pdf}%
\end{picture}%
\setlength{\unitlength}{3947sp}%
\begingroup\makeatletter\ifx\SetFigFont\undefined%
\gdef\SetFigFont#1#2#3#4#5{%
  \reset@font\fontsize{#1}{#2pt}%
  \fontfamily{#3}\fontseries{#4}\fontshape{#5}%
  \selectfont}%
\fi\endgroup%
\begin{picture}(2111,1864)(4501,-5519)
\put(4501,-4561){\makebox(0,0)[lb]{\smash{{\SetFigFont{12}{14.4}{\rmdefault}{\mddefault}{\updefault}{\color[rgb]{0,0,0}$v_2$}%
}}}}
\put(4951,-4561){\makebox(0,0)[lb]{\smash{{\SetFigFont{12}{14.4}{\rmdefault}{\mddefault}{\updefault}{\color[rgb]{0,0,0}$w_2$}%
}}}}
\put(5326,-4261){\makebox(0,0)[lb]{\smash{{\SetFigFont{12}{14.4}{\rmdefault}{\mddefault}{\updefault}{\color[rgb]{0,0,0}$v_1$}%
}}}}
\put(5326,-3811){\makebox(0,0)[lb]{\smash{{\SetFigFont{12}{14.4}{\rmdefault}{\mddefault}{\updefault}{\color[rgb]{0,0,0}$w_1$}%
}}}}
\put(5251,-5011){\makebox(0,0)[lb]{\smash{{\SetFigFont{12}{14.4}{\rmdefault}{\mddefault}{\updefault}{\color[rgb]{0,0,0}$v_1$}%
}}}}
\put(5251,-5461){\makebox(0,0)[lb]{\smash{{\SetFigFont{12}{14.4}{\rmdefault}{\mddefault}{\updefault}{\color[rgb]{0,0,0}$w_1$}%
}}}}
\put(6151,-4561){\makebox(0,0)[lb]{\smash{{\SetFigFont{12}{14.4}{\rmdefault}{\mddefault}{\updefault}{\color[rgb]{0,0,0}$v_2$}%
}}}}
\put(5701,-4561){\makebox(0,0)[lb]{\smash{{\SetFigFont{12}{14.4}{\rmdefault}{\mddefault}{\updefault}{\color[rgb]{0,0,0}$w_2$}%
}}}}
\put(4651,-3886){\makebox(0,0)[lb]{\smash{{\SetFigFont{12}{14.4}{\rmdefault}{\mddefault}{\updefault}{\color[rgb]{0,0,0}$0$}%
}}}}
\put(6076,-5386){\makebox(0,0)[lb]{\smash{{\SetFigFont{12}{14.4}{\rmdefault}{\mddefault}{\updefault}{\color[rgb]{0,0,0}$c$}%
}}}}
\put(6001,-3811){\makebox(0,0)[lb]{\smash{{\SetFigFont{12}{14.4}{\rmdefault}{\mddefault}{\updefault}{\color[rgb]{0,0,0}$a$}%
}}}}
\put(4651,-5386){\makebox(0,0)[lb]{\smash{{\SetFigFont{12}{14.4}{\rmdefault}{\mddefault}{\updefault}{\color[rgb]{0,0,0}$b$}%
}}}}
\end{picture}%

%% file: pyramid_checkerboard.pdf_t
\begin{picture}(0,0)%
\includegraphics{pyramid_checkerboard.pdf}%
\end{picture}%
\setlength{\unitlength}{3947sp}%
\begingroup\makeatletter\ifx\SetFigFont\undefined%
\gdef\SetFigFont#1#2#3#4#5{%
  \reset@font\fontsize{#1}{#2pt}%
  \fontfamily{#3}\fontseries{#4}\fontshape{#5}%
  \selectfont}%
\fi\endgroup%
\begin{picture}(3786,3828)(7835,-175)
\put(7843,-15){\makebox(0,0)[lb]{\smash{{\SetFigFont{7}{8.4}{\rmdefault}{\mddefault}{\updefault}{\color[rgb]{0,0,0}$q_0,q_c$}%
}}}}
\put(7835,473){\makebox(0,0)[lb]{\smash{{\SetFigFont{7}{8.4}{\rmdefault}{\mddefault}{\updefault}{\color[rgb]{0,0,0}$q_b,q_a$}%
}}}}
\put(8401,-112){\makebox(0,0)[lb]{\smash{{\SetFigFont{7}{8.4}{\rmdefault}{\mddefault}{\updefault}{\color[rgb]{0,0,0}$q_a,q_b$}%
}}}}
\end{picture}%